\documentclass[10pt]{article}

\usepackage{amsmath}
\usepackage{amssymb}

\usepackage{graphicx}

\usepackage{cite}
\usepackage{subfigure}
\usepackage{color} 

\usepackage{float} 

\usepackage[labelfont=bf,labelsep=period,justification=raggedright]{caption}

\usepackage{hyperref}

\DeclareMathOperator*{\argmax}{argmax}

\makeatletter
\renewcommand{\@biblabel}[1]{\quad#1.}
\makeatother

\date{2010-08-27}

\begin{document}

\begin{flushleft}
{\Large
\textbf{Optimizing topological cascade resilience based on the structure of terrorist networks}
}
\\
Alexander Gutfraind$^{1}$, 
\\
\bf{1} Center for Nonlinear Studies and T-5/D-6, Los Alamos National Laboratory, Los Alamos, New Mexico, USA 87545
\\
$\ast$ E-mail: agutfraind.research$@$gmail.com
\end{flushleft}

\section{Abstract}
Complex socioeconomic networks such as information, finance
and even terrorist networks need resilience to cascades - to prevent
the failure of a single node from causing a far-reaching domino effect.
We show that terrorist and guerrilla networks are uniquely cascade-resilient
while maintaining high efficiency, but they become more vulnerable
beyond a certain threshold. We also introduce an optimization
method for constructing networks with high passive cascade resilience.
The optimal networks are found to be based on cells, where each cell has a star topology. 
Counterintuitively, we find that there are conditions where networks
should not be modified to stop cascades because doing so would come at
a disproportionate loss of efficiency. Implementation of these findings
can lead to more cascade-resilient networks in many diverse areas.
Keywords: complex networks, resilience, cascade, multi-objective optimization, epidemics on networks, terrorist networks.


\section{Introduction}
Cascades are ubiquitous in complex networks and they have inspired
much research in modeling, prediction and mitigation \cite{Pastor-Satorras01,Crepey06,Centola07,Huang07,Buldyrev10,Newman02,Newman02virus,Motter02,Watts02,Motter04}.
For example, since many infectious diseases spread over contact networks a single carrier
might infect other individuals with whom she interacts.
The infection might then propagate widely through the network, leading to an epidemic.
Even if no lives are lost, recovery may require both prolonged hospitalizations and expensive treatments.
Similar cascade phenomena are found in other domains such as power distribution systems \cite{Dobson07,Lai04,Johnson10},
computer networks such as ad-hoc wireless networks \cite{Newman02virus},
financial markets \cite{Battiston07,Iori08} and socio-economic systems \cite{Kempe03}.
A particularly interesting class are ``dark'' or clandestine social networks,
such as terrorist networks, guerrilla groups \cite{Raab03}, espionage and crime rings \cite{Baker93,Morselli07}.
In such networks if one of the nodes (i.e. individuals) is captured by law enforcement
agencies, he may betray all the nodes connected to him leading to their
likely capture.

\begin{figure*}[!h]
\begin{centering}
\includegraphics[clip,width=1.0\columnwidth]{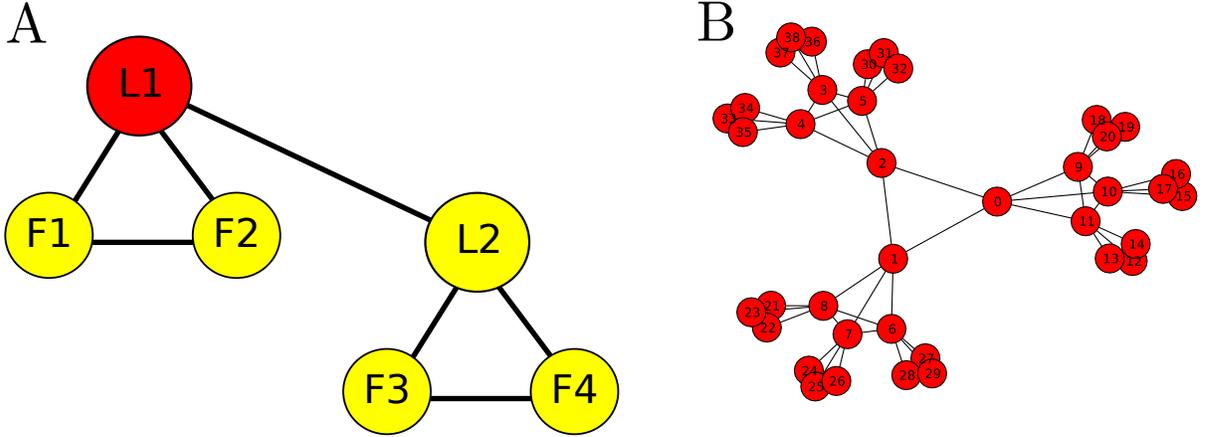}
\par\end{centering}

\caption{{\bf The French World-War II underground network \emph{Francs-tireurs et Partisans }(FTP) reconstructed by the author based on the account in \cite{Miksche50}.} Its organizational unit was the combat group (A). In an idealized case, nor always followed, this was divided into two {}``teams'' of three fighters, where leader L1 was in overall command and in command of team $1$. His lieutenant, L2, led team $2$ and assumed overall command if L1 was captured. The small degree of the nodes ensured that the capture of any one node did not risk the exposure of a significant fraction of the organization. Each {}``group'' is in a command hierarchy (B) where $3$ groups (bottom-level nodes) made a {}``section'', $3$ sections made a {}``company'', and finally $3$ companies made a {}``battalion''.} \label{fig:cellular}
\end{figure*}

Dark networks are therefore designed to operate in conditions of intense cascade pressure.
As such they might serve as useful prototypes
of networks that are cascade-resilient because of their connectivity structure (topology) alone.
Their nodes are often placed in well-defined cells - closely-connected subnetworks with only sparse connections to the outside
(for an example from World War II see Fig.~\ref{fig:cellular}) \cite{Miksche50}.
The advantages of cells are thought to be 
that the risk from the capture of any person is mostly limited to his or her cell mates,
thereby protecting the rest of the network \cite{Gunther78,Lindelauf08hetero}.
Modern terrorist groups retain this cellular structure, but increasingly
use networks made of components with no connections between them, thus 
caging cascades within each component \cite{Rodriguez04,Sageman08,Woo09}.

To represent networks from different domains, this paper will use simple unweighted graphs.
This approach offers simplicity and can employ tools from the well-developed field of graph theory.
A simplification is also unavoidable given the lack of data on networks,
especially on dark networks where only the connectivity is known, if that.
Ultimately through, models of networks, especially dark networks must
consider their evolving nature, fuzzy boundaries and multiplicities
of node classes and diverse relationships.

Fortunately, the loss of information involved in representing networks as simple
rather than as weighted graphs could be evaluated.
In the Supporting Information section we consider two unusually rich data sets 
where the edges could be assigned weights.
We find that the error in using simple graphs has no systematic bias and is usually small.

\subsection{Evaluating Cascade Resilience of Networks}
Our preliminary task is to compare the cascade resilience of networks from different domains.
We will see that dark networks are indeed more successful in the presence of cascades than other complex networks.
Their success stems not from cascade resilience alone
but from balancing resilience with efficiency (a measure of their ability to serve their mission.) 

We will consider a particular type of cascade resilience and a particular definition of efficiency.
For resilience we will use a probabilistic process known as ``SIR'' (susceptible-infected-recovered.)
In SIR any failed (captured) node leads to the failure of each neighboring
node independently with probability $\tau$ \cite{Newman03review}.
Using the SIR model, resilience $R(G)$ could be defined as the average fraction
of the network that does not fail in the cascade. Efficiency
$W(G)$ is also a function of the connectivity structure, and could
be defined based on the distances between all pairs of nodes in the
graph (see the Methods section for exact expressions.) 

Observe that the most cascade-resilient
network is the network with no edges (hence no cascades can propagate),
but it is also the least efficient kind of network. 
It is expected that resilience and efficiency will be in opposition,
requiring trade-offs. Just as disconnected networks are resilient and
inefficient, highly-efficient networks such as densely-connected graphs
are likely to have low resilience (for a historic example see~\cite{Zawodny78}.)

Define the overall ``fitness'', $F(G)$, of
a network by aggregating resilience and efficiency through a weight
parameter $r$: \[
F(G)=rR(G)+(1-r)W(G)\,.\]
The parameter $r$ depends on the application and represents the cost
of restoring the network after a cascade - from light ($r\to0$) to
catastrophic ($r\to1$). 
It is possible to include in fitness other metrics such as construction cost.

We will compare the fitnesses of several complex networks, including communication, infrastructure and scientific networks
to the fitnesses of dark networks.
The class of dark networks will be represented by three networks: the 9/11, 11M and FTP networks.
The 9/11 network links the group of individuals who were directly involved in the September 11, 2001 attacks on New York and Washington, DC~\cite{Krebs02}.
Similarly the 11M network links those responsible for the March 11, 2004 train attacks in Madrid~\cite{Rodriguez04}.
Both 9/11 and 11M were constructed from press reports of the attacks. Edges in those networks connect two individuals who worked with each other in the plots~\cite{Krebs02, Rodriguez04}.
The FTP network is an underground group from World War II (Fig.~\ref{fig:cellular}), whose network was constructed by the author from a historical account~\cite{Miksche50}.

\begin{figure*}[!h]
\begin{centering}
\includegraphics[clip,width=1.0\columnwidth]{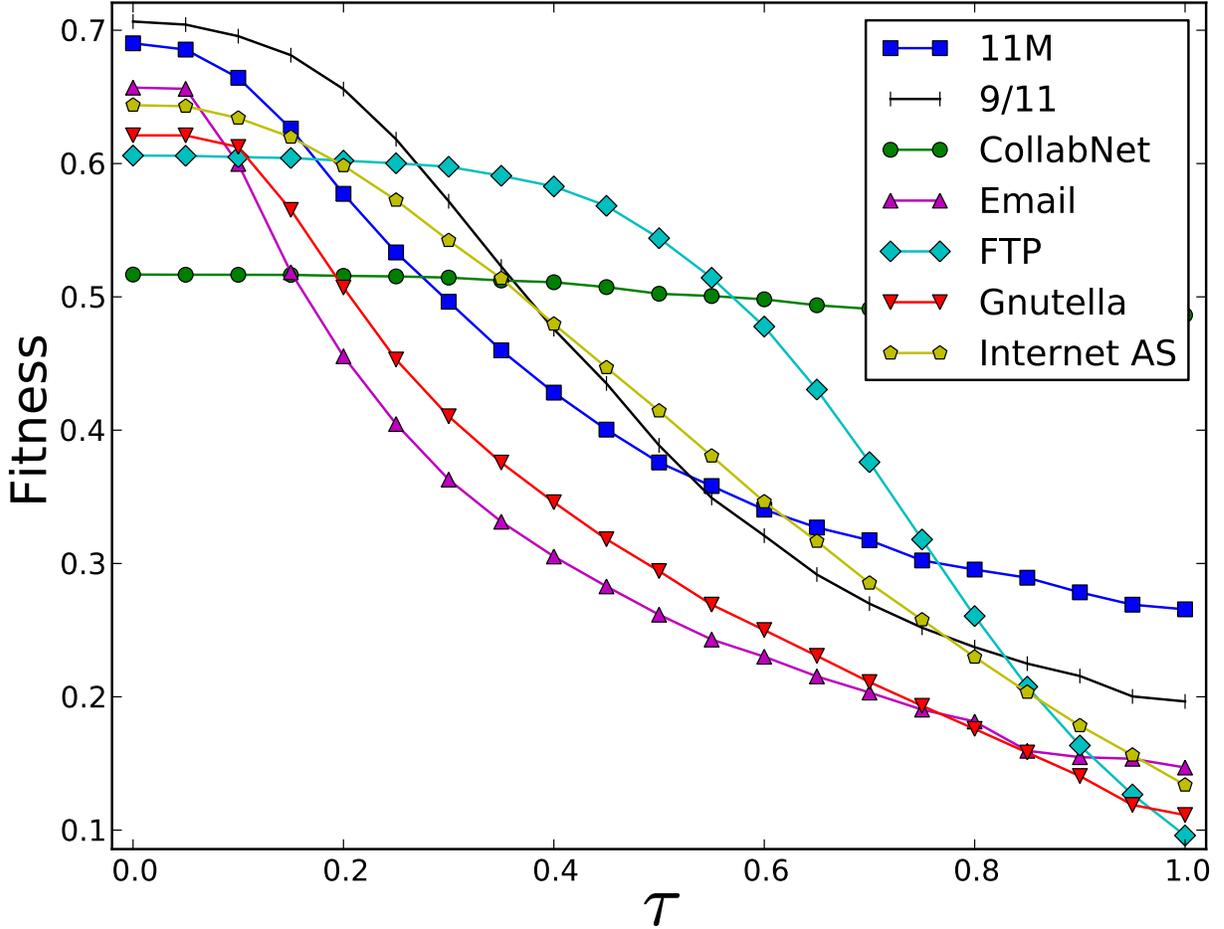}
\par\end{centering}

\caption{{\bf Fitnesses of various networks at $r=0.51$ and various values of $\tau$.}  11M is the network responsible for the March 11, 2004 attacks in Madrid ($70$ nodes, $240$ edges.) 9/11 \cite{Krebs02} is the network responsible for the 9/11 attacks ($62$ nodes, $152$ edges.) CollabNet \cite{Newman04commun} is a scientific co-authorship network in the area of network science ($1589$ nodes, $2742$ edges.) E-Mail \cite{Arenas03email} is a university's e-mail contact network, showing its organizational structure ($1133$ nodes, $5452$ edges.) FTP is the network in Fig.~\ref{fig:cellular} ($174$ nodes, $300$ edges.) Gnutella \cite{Ripeanu02,Leskovec07} is a snapshot of the peer-to-peer network ($6301$ nodes, $20777$ edges.) Internet AS \cite{LeskovecKF05} is a snapshot of the Internet at the autonomous system level ($26475$ nodes, $53381$ edges.) Except for $\tau>0.6$ dark networks (11M, 9/11 and FTP) attain the highest fitness.} 

\label{fig:opt-fitness-emp}
\end{figure*}

Figure~\ref{fig:opt-fitness-emp} shows that the dark networks attain the highest fitness values of all networks, except
for extreme levels of cascade risk ($\tau>0.6$). 
This is to be expected: only 11M, 9/11, and the FTP networks have
been designed with cascade resilience as a significant criterion - a property that makes them useful case studies. 
For high cascade risks ($\tau>0.6$) the CollabNet network exceeds the fitnesses of the dark networks.
CollabNet was drawn by linking scientists who co-authored a paper in the area of network science \cite{Newman04commun}.
It achieved high fitness because it is partitioned into research groups that have no publications with outside scientists. 
Like some terrorist networks, it is separated into entirely disconnected cells.

The 9/11 and the 11M networks are very successful for low values of $\tau$ ($<0.2$),
but then rapidly deteriorate because of a jump in the extent of cascades
- the so-called percolation transition \cite{Draief08}.
Past this threshold, cascades start affecting a large fraction of the network, resilience collapses and the fitness declines rapidly. 
The pattern of onset of failure can be clearly seen in most of the networks. 
For violent secret societies this transition means that the network might
be initially hard to defeat, but there is a point after which efforts
against it start to pay off. Because $\tau$ is representative of the security
environment, the 9/11 network is found to be relatively ill-adapted
to the more stringent security regime implemented after the attacks. Indeed,
it is likely that the 9/11 attacks would have been thwarted under
the current security regime since some of the nodes were captured before the attacks,
but not interrogated in time to discover and apprehend the rest of the network \cite{911report}.
In contrast, the cellular tree hierarchy of the FTP network is more
suitable for an intermediate range of cascade risks.  However, the pair-wise
distances in it are too long to provide high efficiency. Therefore,
its fitness is comparatively poor in the very low and very high
values of $\tau$.

\section{Designing Networks}
The success of dark networks must be due to structural elements of those networks, such as cells. 
If identified, those elements could be used to design more resilient networks
and to upgrade existing ones.
Thus, by learning how dark networks organize, it will be possible to make
networks such as communication systems, financial networks, and others more resilient and efficient.

\begin{figure*}
\begin{centering}
\includegraphics[clip,width=1.0\columnwidth]{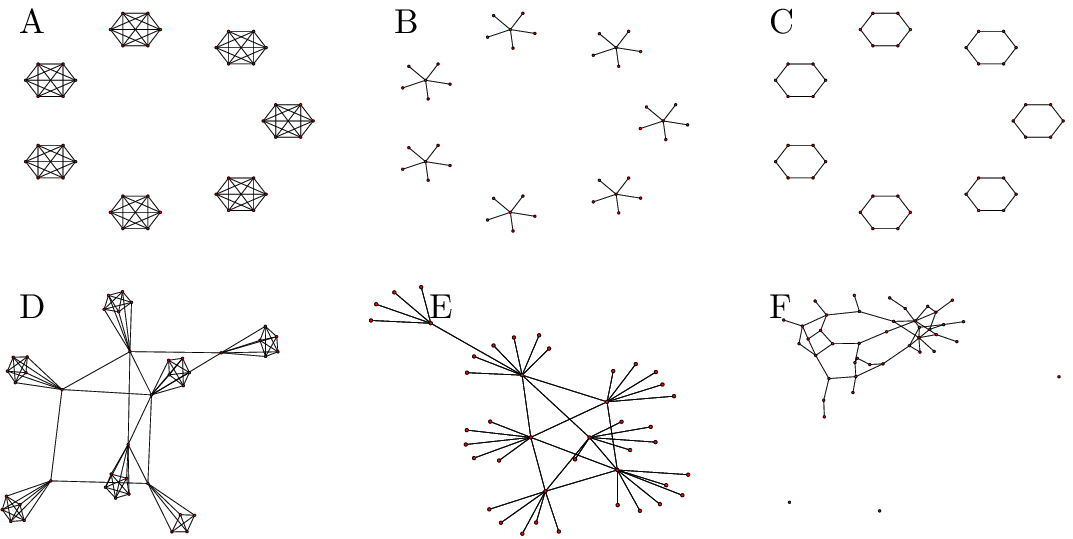}
\par\end{centering}
\caption{{\bf Graphs illustrating the $6$ network designs.} Cliques (A), Stars (B), Cycles (C), Connected Cliques (D), Connected Stars (E), and Erdos-Renyi ``ER'' (F). Each design is configured by just one or two parameters (the number of individuals per cell and/or the random connectivity.) This enables rapid solution of the optimization problem. In computations the networks were larger ($n=180$ nodes.)} \label{fig:designs-1}
\end{figure*}

Those identification and design problems are our next task.
We propose to solve both using an approach based on discrete optimization. 
Let a set of graphs $\mathbb{G}$ be called a {}``network design'' if all the
networks in it share a structural element. Since dark networks are often
based on dense cliques, we consider a design where all the networks
consist of one or multiple cliques.
We consider also designs based on star-like cells, cycle-based cells
and more complex patterns (see Fig.~\ref{fig:designs-1} and SI for the exact set of networks.)

In the first step we will find the most successful network within each design.
Namely, consider an optimization problem where the decision variable
is the topology $G$ of a simple graph taken from a design $\mathbb{G}$.
The objective is the fitness $F(G)$:
\begin{equation}
\max_{G\in\mathbb{G}}F(G)\,.\label{eq:objective}
\end{equation}
In the second step we will compare the fitnesses across designs,
thus identifying the topological feature with the highest fitness (e.g. star vs. clique.)

This optimization problem could be used more broadly:
It introduces a method for designing cascade-resilient networks for applications
such as vital infrastructure networks.
To apply this to a given application, one must make the
design $\mathbb{G}$ the set of all feasible networks in that domain, to the extent possible
by computational constraints.
In the area of terrorist networks, the model is closely related to the game-theoretic
work of Lindelauf et al.\cite{Lindelauf08covert,Lindelauf08hetero}.

A complementary approach is to consider the multi-objective optimization problem
in which $R(G)$ and $W(G)$ are maximized simultaneously:
\begin{equation}
\max_{G\in\mathbb{G}}\left\{ R(G),W(G)\right\} \,.\label{eq:multi-objective}
\end{equation}
The multi-objective approach cannot find the optimal network but instead
produces the Pareto frontier of each design - the set of network configurations
that cannot be improved without sacrificing either efficiency or
resilience. The decision maker can use the frontier to make the optimal
trade-off between resilience and efficiency.

\section{Results}
\subsection{Optimal Network\label{sub:Optimal-Network}}

\begin{figure*}[!h]
\begin{centering}
\includegraphics[width=1\columnwidth]{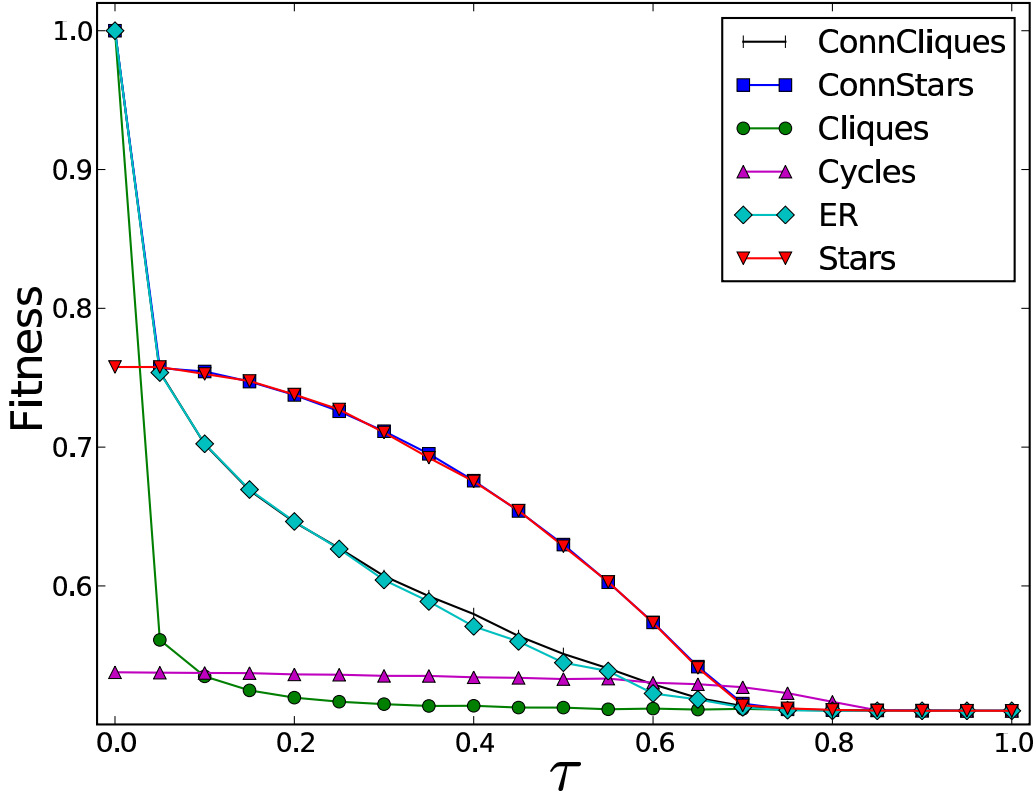}
\par\end{centering}
\caption{{\bf Fitness at $r=0.51$ of various network designs.}
The Connected Stars design is the best design at all cascade risks, $\tau$.
Cliques and Connected Cliques are competitive only for extreme ranges of $\tau$. The superiority of Connected Stars over the ER (random graph) 
confirms the hypothesis that cells give fitness gains against cascades.
The fitness of a design at each value of $\tau$ is defined as the fitness of the optimal configuration (network ensemble) within that design.
}
\label{fig:opt-fitness}
\end{figure*}

The first set of experiments compares the designs against each other
under different cascade risks ($\tau$), Fig.~\ref{fig:opt-fitness}.
At each setting of $\tau$, each design is optimized to its best configuration,
i.e. the best cell size, and connectivity if applicable. The curves
indicate the fitness of the optimal network in each design.  Typically
at each $\tau$ the optimal network is different from the optimal network at another $\tau$.
Observe that within each design, as $\tau$ increases the fitness decreases - one
cannot win when fighting cascades, only delay (see SI for proof.) In certain applications it
is possible to invest in reducing the cascade propagation probability,
$\tau$. Then the curves in Fig.~\ref{fig:opt-fitness} could also
be viewed as expressing the gain from efforts to reduce cascades by
reducing $\tau$ and also adapting the network structure. If the slope
is steep then the gains are large. 

Comparing designs to each other reveals that Connected Stars is
superior to all others in fitness (Fig.~\ref{fig:opt-fitness}.)
The design also outperforms any of the empirical networks in Fig.~\ref{fig:opt-fitness-emp}
in part because for each value of $\tau$ we selected the optimal network.
The simpler Stars design is almost as fit, deteriorating only at extreme
ranges of $\tau$. The rankings of the designs are of course dependent
on the parameter values, but not strongly (see SI for proof.) 
Star-like designs are successful
because the central node in a star acts as a cascade blocker while
keeping the average distance in the star short ($\sim2$). Only for
sufficiently low $r$, the Cliques, Connected Cliques and Connected
Stars designs are superior to the Stars design. For such values of
$r$ efficiency is the dominant contributor to fitness.
High weighting for efficiency benefits the former designs where efficiency can be $1$ by constructing a fully
connected (complete) graph (see SI for analytic results.) 
In the star design efficiency is lower, reaching $\sim\frac{1}{2}$ (when all nodes are placed in a single large star.)

It has been long conjectured that cells provide dark networks with high resilience.
Indeed, this is probably the reason why we found that dark networks have higher fitnesses than other networks.
But cells also reduce the efficiency of a network since they isolate nodes from each other.
To rigorously determine the net effect of cells, we compare the ER design (random graphs) to the Connected Stars design.
ER is a strict subset of Connected Stars but only Connected Stars has cells.
Therefore it is notable that Connected Stars has a higher fitness than ER, often significantly so.
Indeed, cells must be the cause of higher fitness because cells are the only feature in Connected Stars that ER lacks.

\subsection{Properties of Optimal Networks}

\begin{figure*}[h]
\begin{centering}
\includegraphics[clip,width=1.0\columnwidth]{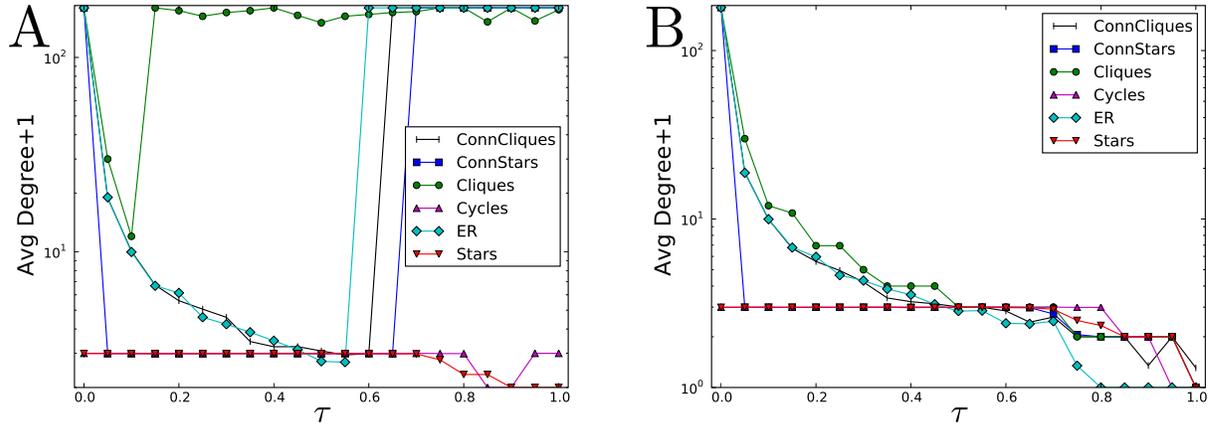}
\par\end{centering}
\caption{{\bf Average degree in the optimal configuration of each design.} At $r=0.49$ (A) the optimization prefers networks that have
high efficiency while at $r=0.51$ (B) the preference is for resilience.  In (B) the average degree diminishes monotonically
to compensate for increasing cascade risk. In (A) most designs have a threshold $\tau$ at which they jump back to a completely-connected graph
because structural cascade resilience becomes too expensive in terms of efficiency.
}
\label{fig:opt-avgDeg}
\end{figure*}

Many properties of the optimal networks such as resilience, efficiency
and edge density show rapid phase transitions as $r$ is changed.
For example, in the Cliques design when $r<0.5$ the optimal network
has high density that maximizes efficiency, whereas for $r>0.5$ it
is sparse and maximizes resilience (Fig.~\ref{fig:opt-avgDeg}.)

Intuition may suggest that the networks grow more sparse as cascade
risk grows. Instead, the trend was non-monotonic (Fig.~\ref{fig:opt-avgDeg}.)
For $\tau\gg0$ and $r<0.5$ Cliques, Connected Cliques and Connected
Stars became denser, instead of sparser, and for them the most sparse
networks were formed in the intermediate values of $\tau$ where the
optimal networks achieve both relatively high resilience and high
efficiency. At higher $\tau$ values, when $r<0.5$ it pays to sacrifice
resilience because fitness is increased when efficiency is made larger
through an equal or lesser sacrifice in resilience. The Stars design
does not show a transition at $r=0.5$ because it
is hard to increase efficiency with this design.

\subsection{Multi-objective Optimization\label{sub:Multi-objective-Optimization}}

\begin{figure*}[!h]
\begin{centering}
\includegraphics[width=1\columnwidth]{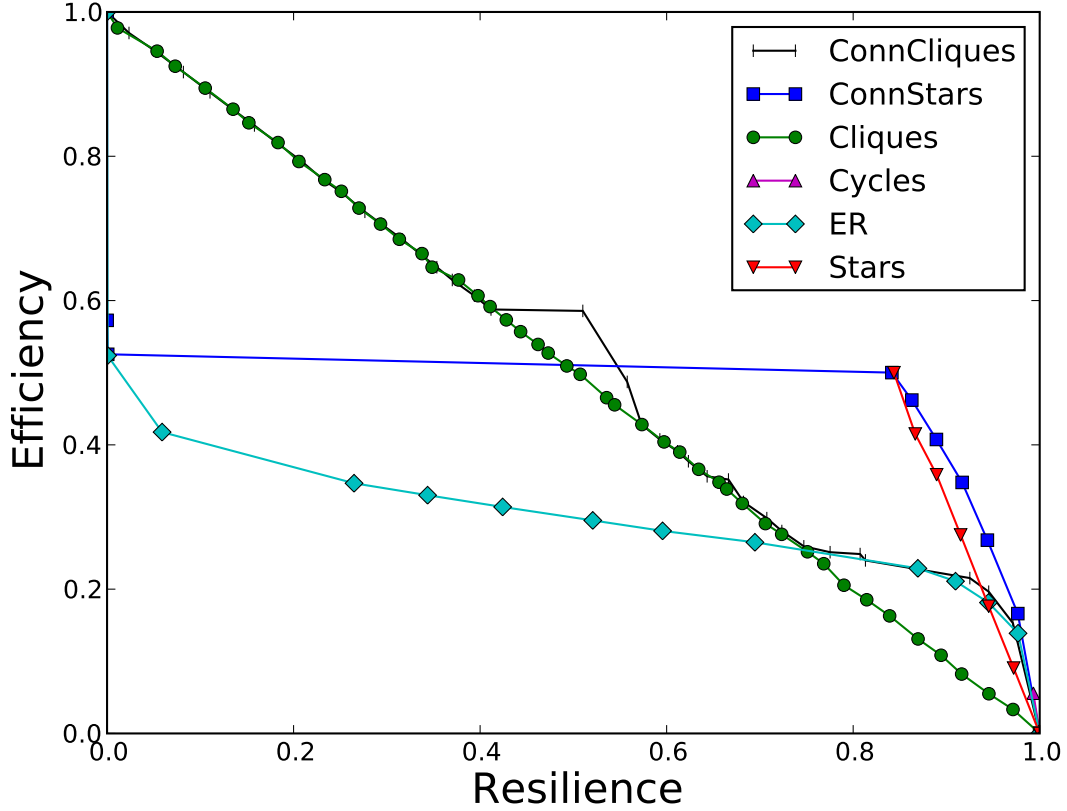}
\par\end{centering}

\caption{{\bf The Pareto frontiers of various network designs ($\tau=0.4$).} The configurations of the Connected Stars design dominate
over other designs when the network must achieve high resilience.  However, designs based on cliques are dominant when
high efficiency is required.  Several designs show sharp transitions where at a small sacrifice of efficiency it is possible
to achieve large increases in cascade resilience.
\label{fig:Pareto-frontier}}

\end{figure*}

A complementary perspective on each design is found from its Pareto frontier
of resilience and efficiency (Fig.~\ref{fig:Pareto-frontier}.) 
Typically a design is dominant in a part of the Resilience-Efficiency plane but not all of it.
The Stars and Connected Stars designs can access most of the high resilience-low efficiency region.
In contrast, the Cliques and Connected Cliques can make networks in the medium resilience-high efficiency regions.

The sharp phase transitions discussed earlier are seen clearly: along most of the frontiers, if
we trace a point while decreasing resilience, there is a threshold
at which a small sacrifice in resilience gives a major gain of efficiency.
More generally, consider the points where the frontier is smooth.
By taking two nearby networks on the frontier one can define 
a rate of change of efficiency with respect to resilience: $\left|\frac{\Delta W}{\Delta R}\right|$.
The ratio can be used to optimize the network without using the parameter
$r$. When $\left|\frac{\Delta W}{\Delta R}\right|\gg1$ the network
optimizer should choose to reduce to the resilience of the network
in order to achieve great gains in efficiency; when $\left|\frac{\Delta W}{\Delta R}\right|\ll1$
efficiency should be sacrificed to improve resilience.

\section{Discussion}

The analysis above considered both empirical networks and synthetic ones.
The latter were constructed to achieve structural cascade resilience and efficiency.
In contrast, in many empirical networks the structure emerges through an unplanned growth process
or results from optimization to factors such as cost rather than blocking cascades.
Without exception the synthetic networks showed higher fitness values despite the fact
that they were based on very simple designs.
This suggests that network optimization can significantly improve the fitness and cascade resilience of networks. 
It means that such an optimization process can indeed be an effective method
for designing a variety of networks and for protecting existing networks from cascades.

Many empirical networks also have power-law degree distributions \cite{Newman03review}.
Unfortunately, this feature significantly diminishes their cascade resilience: 
the resulting high-degree hubs make the networks
extremely vulnerable to cascades once $\tau$ is slightly larger than
$0$ \cite{Pastor-Satorras01,Crepey06}.

In some successful synthetic networks the density of edges increased when the 
cascade risk $\tau$ was high.  
This phenomenon has interesting parallels in non-violent social movements
which are often organized openly rather than as secret underground cells 
even under conditions of severe state repression \cite{Sharp03}.
This openness greatly facilitates recruitment and advocacy, justifying
the additional risk to the participants, just like the sacrifice of
resilience to gain higher efficiency is justified under $r<0.5$ conditions. 

There are other important applications of this work, such as the design of power distribution systems.
For power networks, the definition of resilience and efficiency will need to be changed.
It would also be necessary to use much broader designs and optimization under design constraints such as cost.
Furthermore, this work could also be adapted to domains of increasing concern such as financial
credit networks, whose structure may make them vulnerable to bankruptcies~\cite{Battiston07,Iori08}.

\section{Methods}
\subsection{Measuring Resilience}
Research on graph theory has led to the development of a variety of
metrics of robustness or resilience \cite{Klau05} but here unlike
in many other studies the interest is in resilience to cascades and
not to disconnection. One particularly important and well-characterized
class of cascades are those that start at a single node and then spread
probabilistically to neighboring nodes possibly reaching a large fraction
of the network, termed the SIR model and percolation\cite{Newman03review}. Under
this model, resilience can be defined based on the expected size of
the surviving network: \begin{equation}
R(G)=1-\frac{1}{n-1}\mathbb{E}[\mbox{extent of a cascade}]\,,\label{eq:resilience}\end{equation}
where {}``extent of a cascade'' refers to the ultimate number of
new cases created by a single failed node. For simplicity, cascades
are assumed to start at all nodes with uniform probability.

\subsection{Measuring Efficiency}
For many applications the distance between pairs of nodes in the network
is one of the most important determinants of the network's efficiency
(see e.g. \cite{Latora01,Motter02range,Lindelauf08covert}.) When
nodes are separated by short distances they can easily communicate
and distribute resources to each other. This idea motivates the following
{}``distance-attenuated reach'' metric. For all pairs of nodes $u,v\in V$,
weigh each pair by the inverse of its internal distance (the number
of edges in the shortest path from $u$ to $v$) taken to power $g$:
\begin{equation}
W(G)=\frac{1}{n(n-1)}\sum_{u\in V}\sum_{v\in V\smallsetminus\{u\}}\frac{1}{d(u,v)^{g}},\label{eq:efficiency}
\end{equation}
Normalization by $n(n-1)$ ensures that $0\leq W(G)\leq1$, and only
the complete graph achieves $1$. As usual, for any node $v$ with
no path to $u$, set $\frac{1}{d(u,v)^{g}}=0$. The parameter $g$,
{}``connectivity attenuation{}`` represents the rate at which distance
decreases the connectivity between nodes. In the experiments above
$g=1$.  

Supporting Information below contains detailed information
about the optimization methodology, the simulation process, and sensitivity
as well as rigorous justification of certain claims.

\pagebreak[4]

\appendix
\begin{center}
{\bf \huge Supporting Information}
\end{center}

\section{Resilience and Efficiency of Empirical Networks}

\begin{figure}[H]
\begin{centering}
\includegraphics[width=0.49\columnwidth]{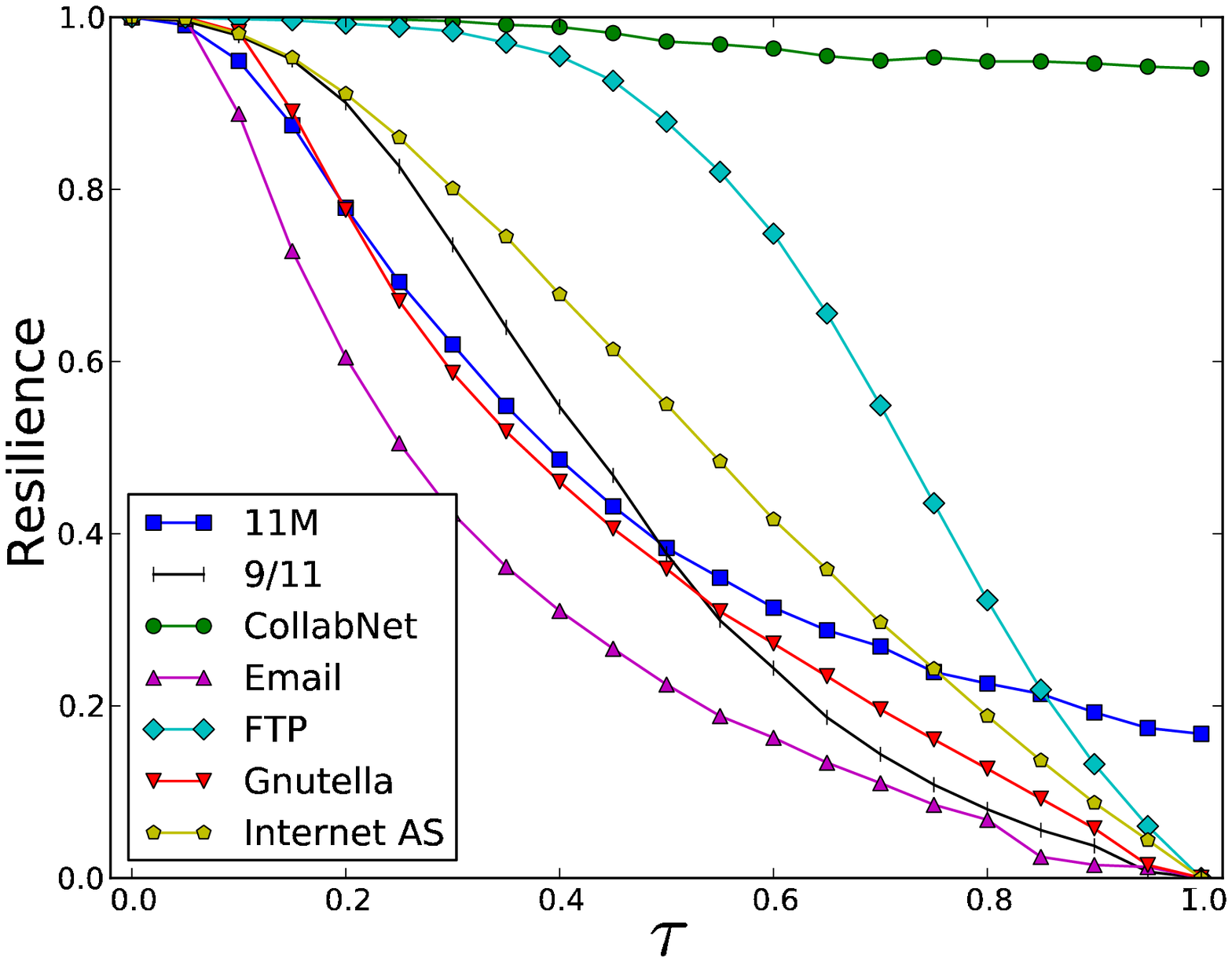}
\includegraphics[width=0.49\columnwidth]{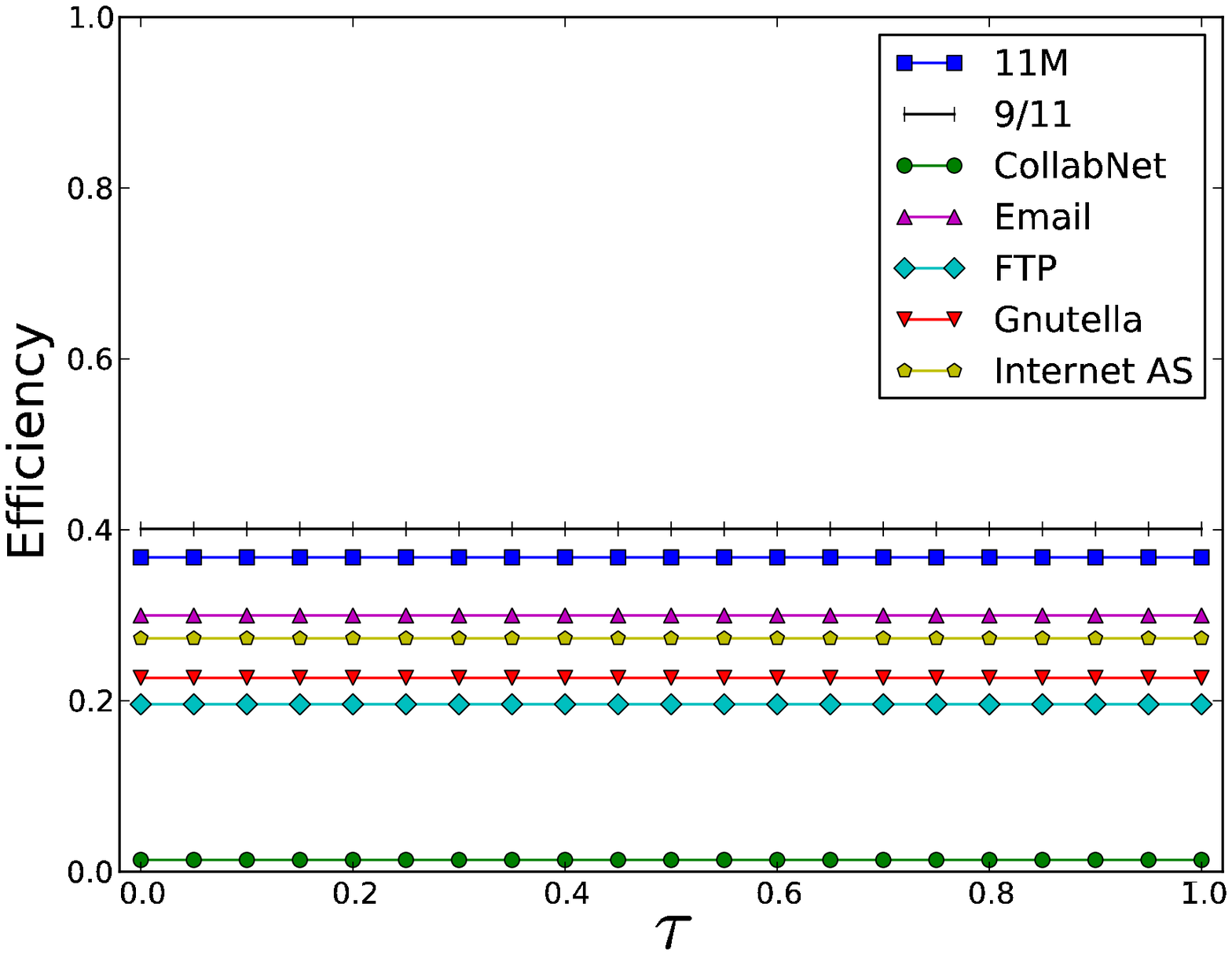}
\par\end{centering}
\caption{{\bf Resilience and Efficiency of the real networks.}  The fittest networks are not always the most resilient.}
\label{fig:additional-emp}
\end{figure}

It is interesting to compare the empirical networks to each other
in their efficiency and resilience (Fig.~\ref{fig:additional-emp}).
Note that FTP and 9/11 networks are not the most resilient, but they
strike a good balance between resilience and efficiency. The advantages
of the two networks over other networks are not marginal, implying
that their advantages in fitness are not sensitive to the choice of
$r$. Of course, they are optimized for particular combinations of
$r$ and $\tau$, and will no longer be very successful outside that
range. For instance, in the range of high $r$ and high $\tau$ networks
with multiple connected components would have higher fitness because
they are able to isolate cascades in one component.

\section{Resilience and Efficiency of Weighted Networks}

In some networks, each edge $(u,v)$ carries a distance weight $D_{uv} > 0$.  
The smaller the distance, the closer the connection between $u$ and $v$.  
To compute the fitness of those networks, we now introduce generalizations of resilience and efficiency.
Those reduce to the original definitions for unweighted networks when $D_{uv}=1$, 
while capturing the effects of weights in the weighted networks.

The original definition of resilience was built on a percolation model where the failure of any
node leads to the failure of its neighbor with probability $\tau$.
In the weighted network more distant nodes should be less likely to spread the cascade.
Thus we make the probability of cascade through $(u,v)$ to be $\min(\tau/D_{uv},1)$.

The efficiency was originally defined as the sum of all-pairs inverse geodesic distances, 
normalized by the efficiency of the complete graph.  In the weighted network, both the distance
and the normalization must be generalized.  To compute the distance $d(u,v)$ we consider the weights
on the edges $D$ and apply Dijkstra's algorithm to find the shortest path.
Normalization too must consider $D$ because a weighted graph with sufficiently small distances
could outperform the complete graph (if all the edges of the latter have $D_{ij}=1$.)
Therefore, we weigh the efficiency by the harmonic mean $H$ of the edges ($E$) of the graph:
\begin{equation}
W(G)=\frac{H(G)}{n(n-1)}\sum_{u\in V}\sum_{v\in V\smallsetminus\{u\}}\frac{1}{d(u,v)^{g}},\label{eq:efficiencyGeneral}
\end{equation}
where $$H(G) = \frac{|E|}{\sum_{(u,v)\in E}{\left(\frac{1}{D_{uv}}\right)^{g}}}\,.$$
The harmonic mean ensures that for any $D$, the complete graph has $W(G)=1$.

Having defined generalized resilience and efficiency we can evaluate the standard
approach to dark network, which represents them as binary graphs $D_{uv} \in \{0,1\}$, rather
than as weighted graphs.  The former approach is often taken because 
the information about dark networks is limited and insufficient to estimate edge weights.

Fortunately, in two cases, the 9/11 network and the 11M network \cite{Krebs02, Rodriguez04} the weights could be estimated.  
The 9/11 data labels nodes as either facilitators or hijackers. 
Hijackers must train together and thus should tend to have a closer relationship.  
Thus set $D_{uv}=2,1,0.5$ if the pair $u,v$ includes zero, one or two hijackers, respectively.
The 11M network is already weighted ($Z_{uv}=1,2,3\dots$) based on the number of functions each contact ($u,v$) 
serves (friendship, kin, joint training etc.) We mapped those weights to $D$ by $D_{uv}=2/Z_{uv}$.
In both networks the transformation was so that the weakest ties have weight $2$, giving them
greater distance than in the binary network, while the strongest ties are shorter than in the binary network.

\begin{figure}[H]
\begin{centering}
\includegraphics[width=0.99\columnwidth]{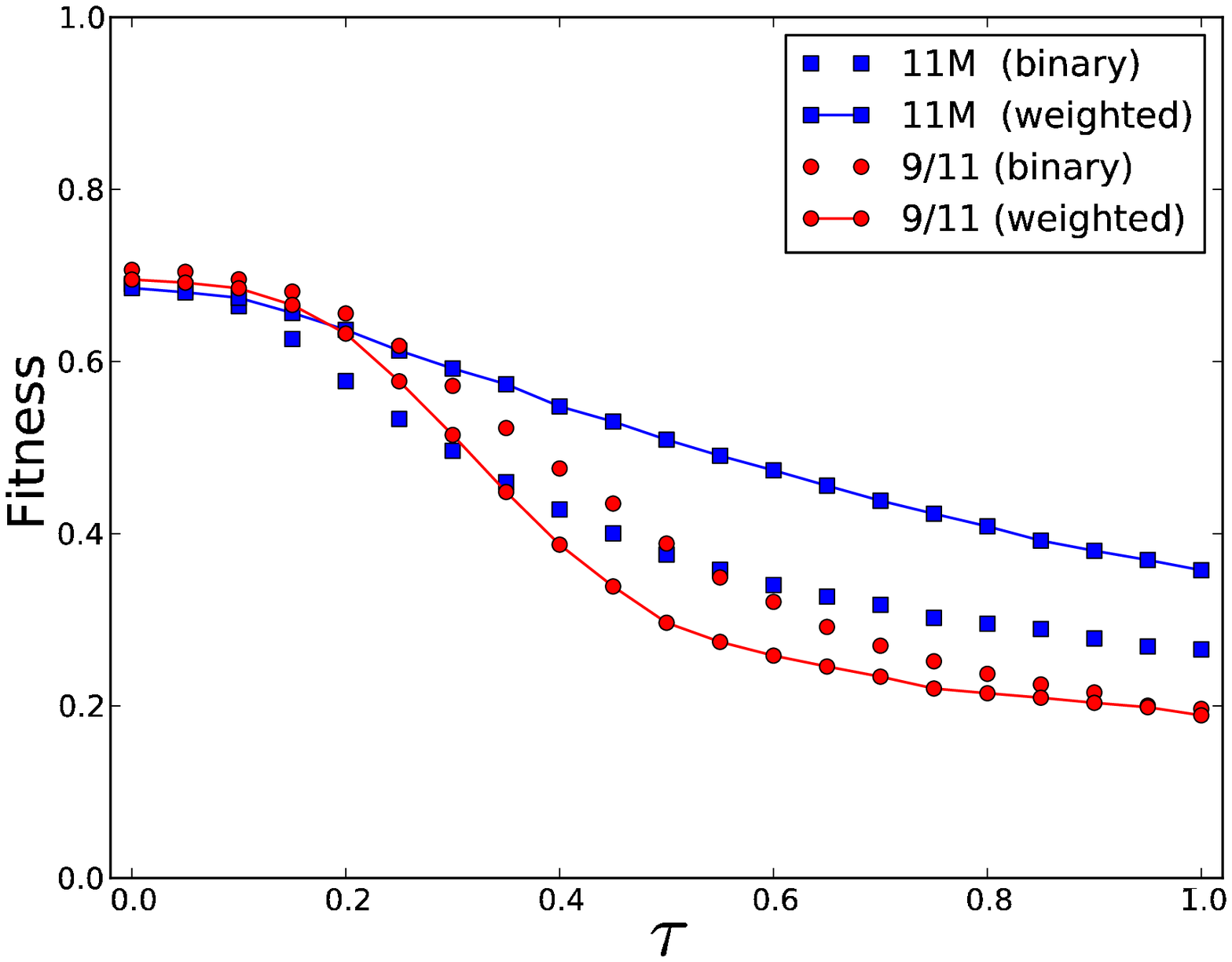}
\includegraphics[width=0.49\columnwidth]{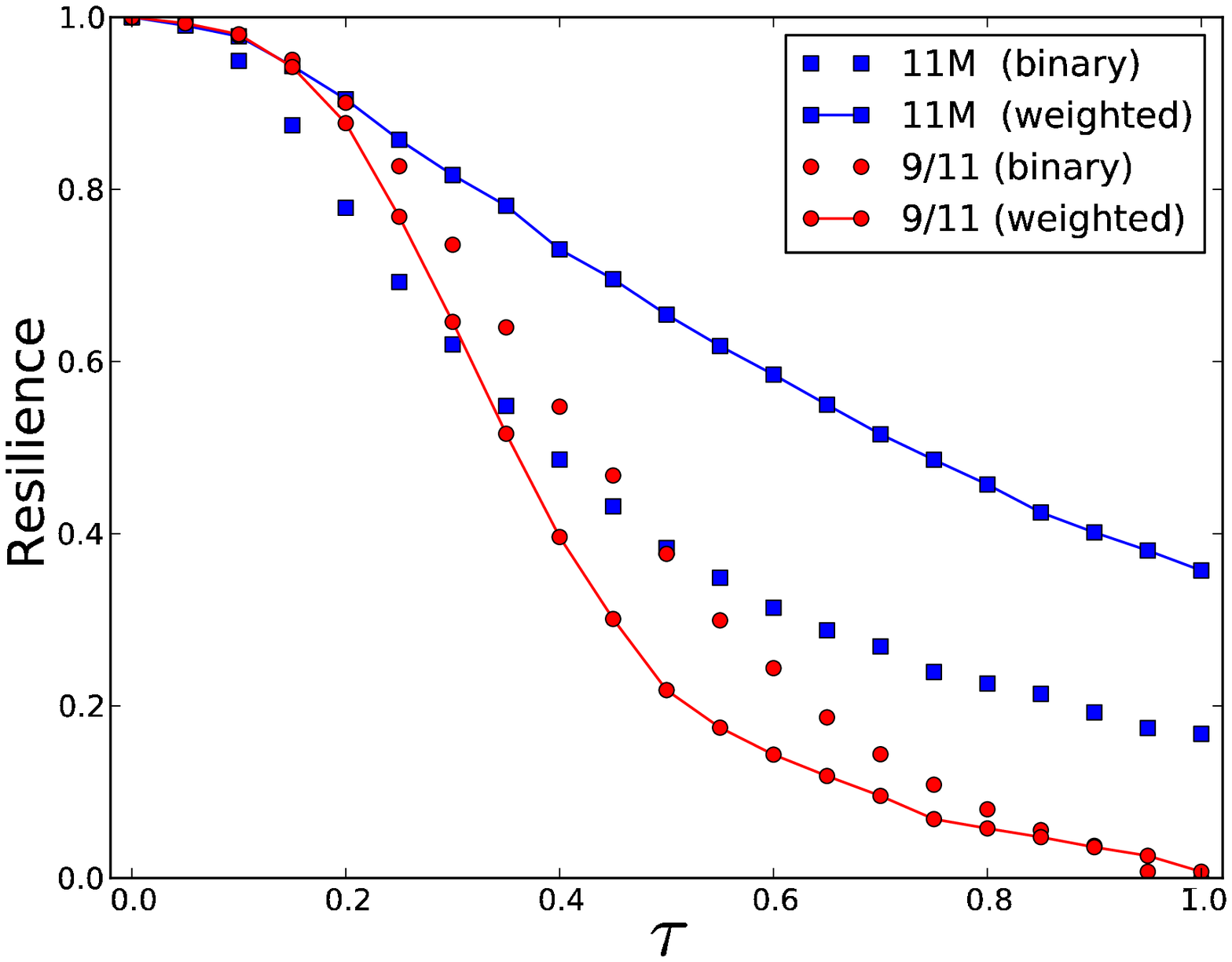}
\includegraphics[width=0.49\columnwidth]{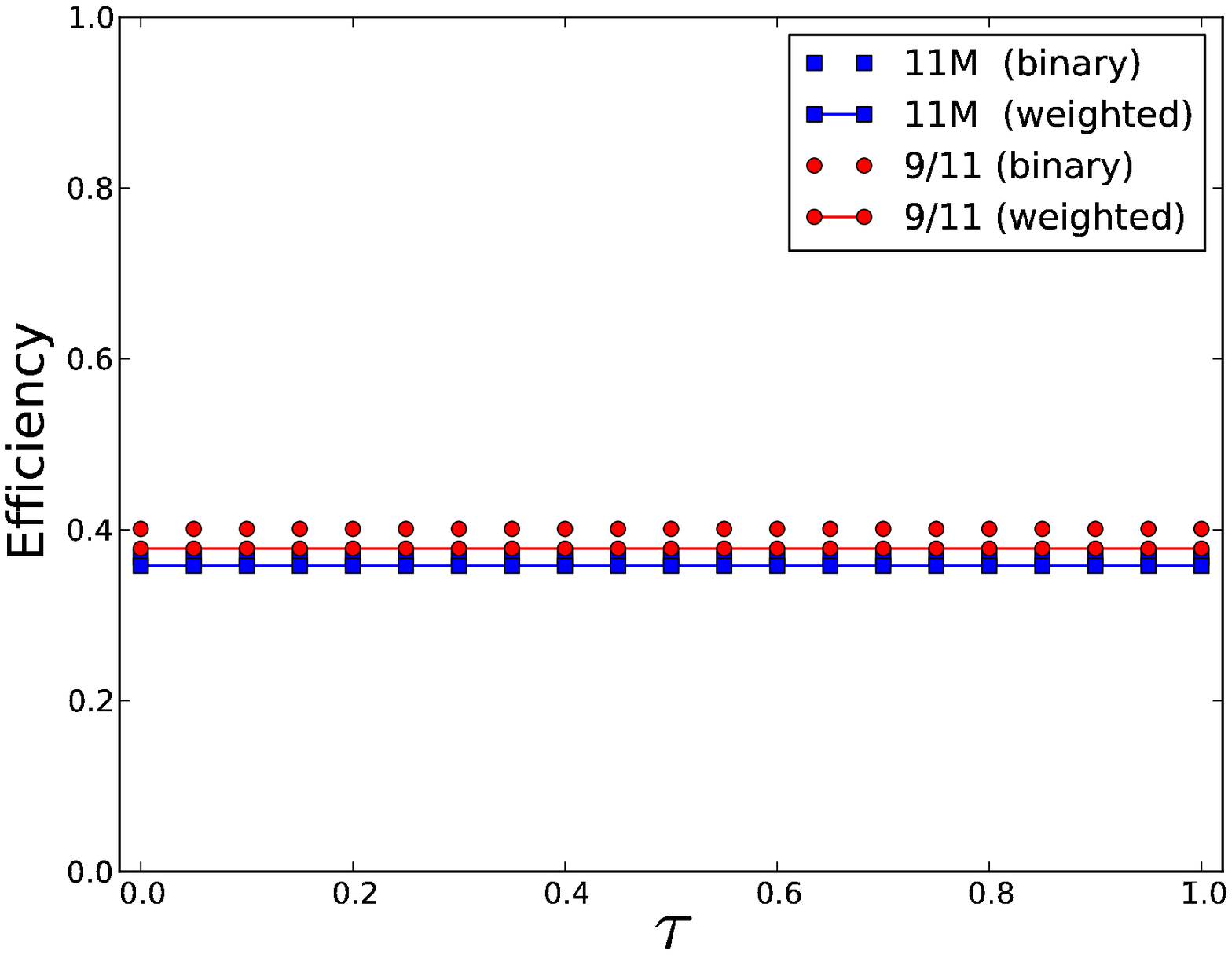}
\par\end{centering}
\caption{{\bf Fitness, Resilience and Efficiency of two dark networks ($r=0.51$), comparing binary and weighted representations.}
The binary representation matches the weighted representation within $0.15$, and typically closer.
\label{fig:binary-emp}
}
\end{figure}

Figure~\ref{fig:binary-emp} compares the fitnesses, resiliences and efficiencies of the weighted and binary representations.
It shows that for both networks, the fitnesses of the binary representation lies within $0.15$ of the fitness of the weighted representation and for some $\tau$ much closer.
The efficiency measures are even more close (within $0.05$.)
The behavior of resilience is intriguing: for the 9/11 network the weighted representation shows
more gradual decline as a function of cascade risk when compared to the binary representation.
For the 11M network the decline is actually slightly more sharp in the weighted representation.
Structurally, the 11M network has a center (measured by betweenness centrality) of tightly knit-nodes (very short distances), 
while the 9/11 network is more sparse at its center, increasing its cascade resilience.  
This effect explains the direction of the error in the binary representation.
Based on those two examples, it appears that the binary representation does not have a systematic bias,
and may even underestimate the fitness of dark networks.

\section{Network Designs\label{sub:Network-Optimization}}

We considered networks on $n=180$ nodes constructed through $6$
simple designs, chosen both based on empirical findings (see e.g.\cite{Arquilla01,Carley06})
as well as the possibility of analytic tractability in some cases.
When more data becomes available on dark networks, it will become
possible to extract additional subgraphs with statistical validity.

Three of the designs are based on identical {}``cells'': each cell
is either (a) a clique (a complete graph), (b) a star (with a central
node called {}``leader''), and (c) a cycle (nodes connected in a
ring). Each of these have a single parameter, $k$ - the number of
nodes in the cell. Recent research suggests that under certain assumptions
constructing networks from identical cells is optimal \cite{Goyal10}.
Let us also consider $n$-node graphs consisting of (d) randomly-connected
cliques (sometimes termed {}``cavemen''), and (e) randomly-connected
stars, in both cases according to probability $p$ . Consider also
(f) the simpler and well-studied Erdos-Renyi (ER) random graph with
probability $p$ (see figure in main text). By considering different
structures for the cells we determine which of those structures provides
the best performance.

The above palette of designs is sufficient for the current study's
purpose of introducing an approach which could be applied to different
domains, as well as begin constructing a theory to address cascade
resilience. There are surely applications where if the objective is
to design cascade-resilient networks, one ought to both reject some
of the designs above for application-specific reasons and also to
introduce additional designs. The task of taking other designs can
be done by merely changing the computer program for generating the
networks. Fruitful future research in cascade-resilient network is
to attack the most general problem of finding the best network on
$n$ nodes. Unfortunately, this search space is exponential and probably
non-convex.

Research on social networks indicates that resilience and efficiency
might be just two of several design criteria that also include e.g.
{}``information-processing requirements'', that impose additional
constraints on network designs \cite{Baker93}. In the original context
{}``information-processing'' refers to the need to have ties between
individuals involved in a particular task, when the task has high
complexity. Each individual might have a unique set of expertise into
which all the other agents must tap directly. Generalizing from sociology,
such {}``functional constraints'' might considerably limit the flexibility
in constructing resilient and efficient networks. Such functional
constraints could be addressed by looking at a palette of network
designs which already incorporate such constraints or using innovative
techniques from convex optimization \cite{Johnson10}. 

The solution to the optimization problem is found by setting each
of the parameters $k$ (and when possible $p$) to various values.
Each design $D$ has {}``configurations'' $C_{1}^{D},C_{2}^{D},\dots$
each specifying the values of the parameters. Each configuration $C_{i}^{D}$
is inputted to a program that generates an ensemble of $1-10$ networks,
whose average performance provides an estimate of the fitness of $C_{i}^{D}$.
The number of networks was $10$ for networks with parameter $p$
because there is higher variability between instances. The coefficient
of variation (CV) in the fitness of the sample networks was monitored
to ensure that the average is a reliable measure of performance. Typically
CV was $<0.2$ except near phase transitions of connectivity and percolation.

Optimization was performed using grid search. Alternative methods
(e.g. Nelder-Mead) were considered but grid search was chosen despite
its computational cost because it suffers no convergence problems
even in the presence of noise (present due to variations in topology
and contagion extent), and collects data useful for sensitivity analysis
and multi-objective optimization. The sampling grid was as follows.
In designs consisting of cells of size $k$, cell size was set to
all integer values in $[1,180]$. If $k$ did not divide $180$, a
cell of size $<k$ was added to ensure that the number of nodes in
the graph is $180$. The number of nodes is $180$ because $180$
is a highly-composite number and so it offers many networks of equally-sized
cells. In general, normalization by $n$ in the definitions of resilience
and efficiency ensures that even when the number of nodes is tripled
the effect of network size on fitness is very small for the above
designs (around $\pm0.05$ in numerical experiments). In designs containing
a parameter of connectivity $p$, it was set to all multiples of $0.05$
in $[0,1]$, with some extra points added to better sample phase transitions.
The grid search algorithm results are readily used to compute the
Pareto frontier using the $\epsilon$-balls method \cite{Laumanns02}
($\epsilon=0.01$).

The resilience metric is most easily computed by simulation where
a node is selected at random to be {}``infected'', and the simulation
is run until all nodes are in states $S$ or $R$, and none is in
state $I$. In the simplest version of the SIR cascade model, which
we adopt, each node in the graph can be in one of three states {}``susceptible'',
{}``infected'' and {}``removed'' designated $S,I$ and $R$ respectively
(these names are borrowed from epidemiology). Time is described in
uniform discrete steps. A node in $S$ state at time $t$ stays in
this state, unless a neighbor {}``infects'' the node, causing it
to move to state $I$ at time $t+1$. Specifically, a node in state
$S$ at time $t$ has probability $\tau$ of turning to $I$ state
at time $t+1$ for each adjacent node in state $I$ at time $t$.
Finally, a node in $I$ state at time $t$ always becomes $R$ at
time $t+1$. Once in state $R$, the node remains there for all future
times. It is possible to consider an alternate model where the rate
of transition $I\to R$ takes more than one time step, but adding
this effect would mostly serve to increase the probability of transmission,
which is already parametrized by $\tau$ \cite{Newman02,Noel09}. 

A cascade/contagion that starts at a single node would run for up
to $n$ steps, but usually much fewer since typically $\tau<1$ and/or
the graph is not connected. To achieve good estimate of the average
extent, the procedure was replicated $40$ times, and then continued
as long as necessary to achieve an error of under $\pm0.5$node with
a $95\%$ confidence interval \cite{Law99}. 

An analytic computation of the cascade extent metric was investigated.
It is possible in theory because the contagion is a Markov process
with states in the superset of the set of nodes, $3^{n}$. Unfortunately,
such a state space is impractically large. When $G$ is a tree, then
an analytic expression exists%
\footnote{Specifically, the mean contagion size is $1+\frac{pG'_{0}(1)}{1-pG'_{1}(1)}$,
where $G_{0}(x)$ generates the degree distribution and $G_{1}(x)=\frac{G'_{0}(x)}{G'_{0}(1)}$
generates the probability of arrival to a node \cite{Newman02}.%
}, and it might be feasible when the treewidth is small \cite{Colcombet02,Newman02}.
However, for many graph designs the tree approximation is not suitable.
Another possible approach is to represent the contagion approximately
as a system of differential equations which can be integrated numerically
\cite{Keeling99} . These possibilities were not pursued since the
simulation approach could be applied to all graphs, while the errors
of the analytic approaches are possibly quite large.

\section{Continuity of Fitness in $r$\label{sec:Continuity-of-Fitness}}

We now justify the claim that fitness is continuous in $r$. In fact,
we will prove the stronger property of Lipschitz-continuity. Notice
that the claim is not about the continuity of fitness of a single
configuration as a function of $r$ but rather that:

\emph{Claim:} $f(r)=\max_{G\in\mathbb{G}}F(G,r)$ is Lipschitz-continuous
for $r\in[0,1]$. 

\emph{Proof:} The argument constructs a bound on the change in $f$
in terms of the change in $r$. Consider an optimal configuration
$C_{1}$ of a design for $r=r_{1}$ and let its fitness be $f_{1}=F(C_{1},r_{1})$
(there is slight abuse of notation since $C$ is a configuration,
whose fitness is the average fitness of an ensemble of graphs).\\
Observation 1: consider the fitness of $C_{1}$ at $r=r_{2}$.
Because $C_{1}$ is fixed and the metrics are bounded ($0\leq R\leq1$
and $0\leq W\leq1$), the fitness change is bounded by the change
in $r$: \begin{eqnarray*}
\left|f_{1}-F(C_{1},r_{2})\right| & = & \left|r_{1}R(C_{1})+(1-r_{1})W(C_{1})\right.\\
 &  & \left.-r_{2}R(C_{1})-(1-r_{2})W(C_{1})\right|\\
 & = & \left|(r_{1}-r_{2})R(C_{1})-(r_{1}-r_{2})W(C_{1})\right|\\
 & \leq & \left|r_{1}-r_{2}\right|\,.\end{eqnarray*}
Observation 2: let $C_{2}$ be the optimal configuration for $r=r_{2}$
and let $f_{2}=F(C_{2},r_{2})$. Since $C_{2}$ is optimal for $r=r_{2}$
it satisfies: $f_{2}\geq f(C_{1},r_{2})$, and so $-f_{2}\leq-F(C_{1},r_{2})$.
It follows that $f_{1}-f_{2}\leq f_{1}-F(C_{1},r_{2})$. Take the
absolute value of the right hand side and apply Observation 1 to get
the bound: $f_{1}-f_{2}\leq\left|r_{1}-r_{2}\right|$. \\
Observation 3: applying the argument of Observations 1\&2 but reversing
the roles of $C_{1}$ and $C_{2}$ implies that $f_{2}-f_{1}\leq\left|r_{1}-r_{2}\right|$. 

Observations 2\&3 give $\left|f_{1}-f_{2}\right|\leq\left|r_{1}-r_{2}\right|$,
proving the result.

\section{Additional Performance Results\label{sec:additional-performance}}

\begin{figure}[H]
\begin{centering}
\subfigure[$r=0.25$]{

\includegraphics[width=0.45\columnwidth]{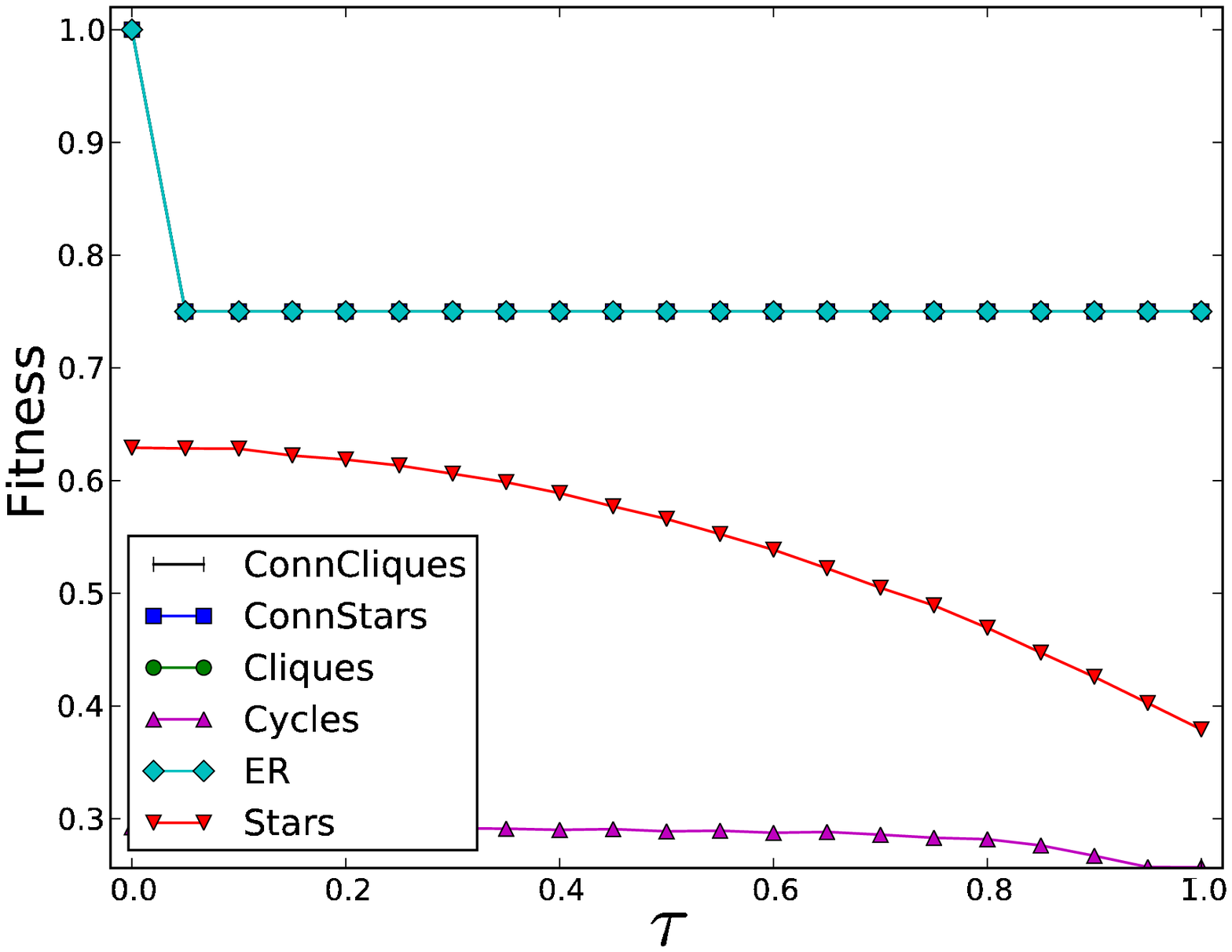}} \subfigure[$r=0.49$]{

\includegraphics[width=0.45\columnwidth]{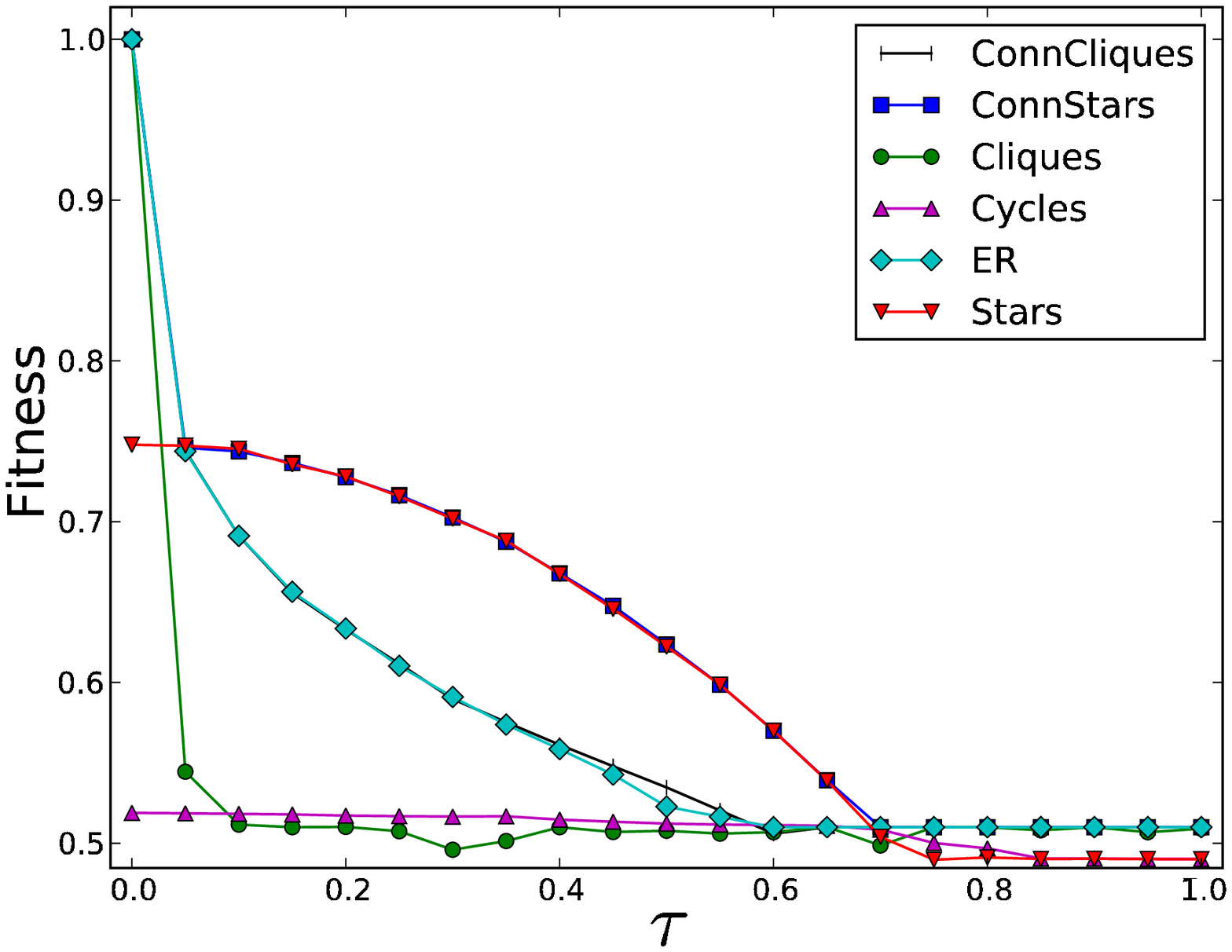}}
\par\end{centering}

\begin{centering}
\subfigure[$r=0.51$]{

\includegraphics[width=0.45\columnwidth]{9_nh_nest_u_gfriend_netstorm_paper_netstorm_hom___e_fitness_R_0d51_2010_04_13__17_21_19_multi}} \subfigure[$r=0.75$]{

\includegraphics[width=0.45\columnwidth]{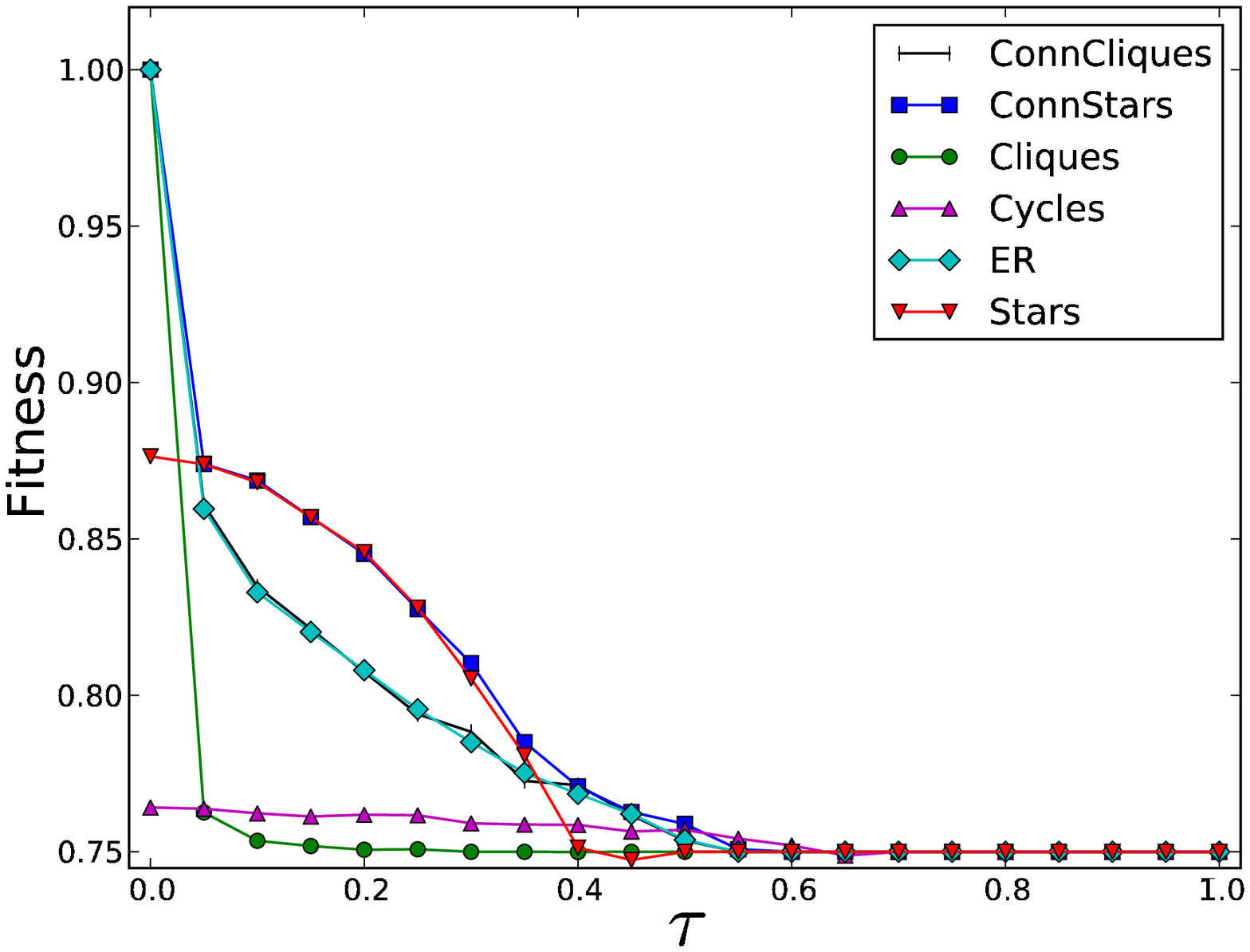}}
\par\end{centering}

\caption{{\bf Fitness of the optimal configuration.} For very low $r$ Cliques outperform Stars.}
\label{fig:fitness-full}
\end{figure}

\begin{figure}[H]
\begin{centering}
\subfigure[$r=0.49$]{

\includegraphics[width=0.45\columnwidth]{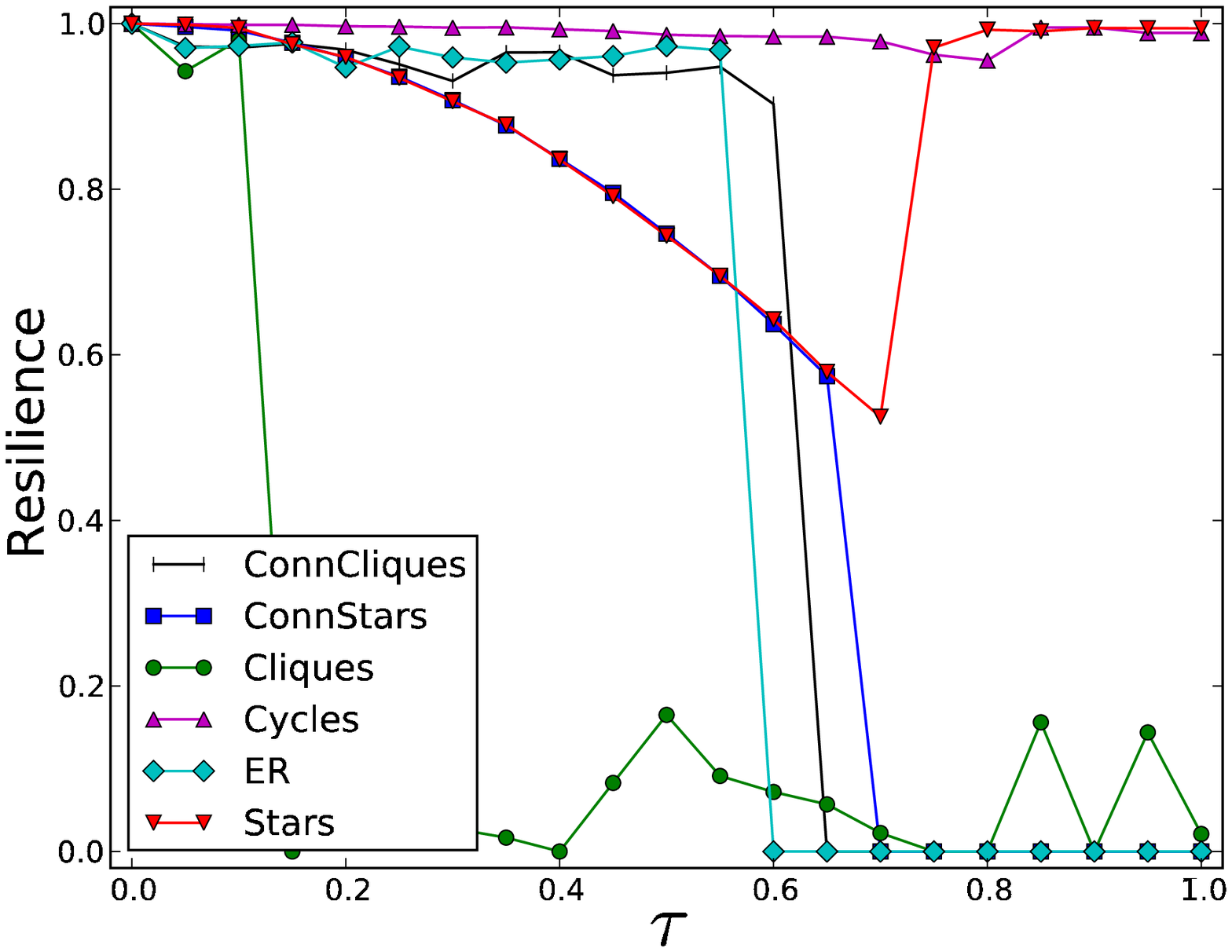}} \subfigure[$r=0.51$]{

\includegraphics[width=0.45\columnwidth]{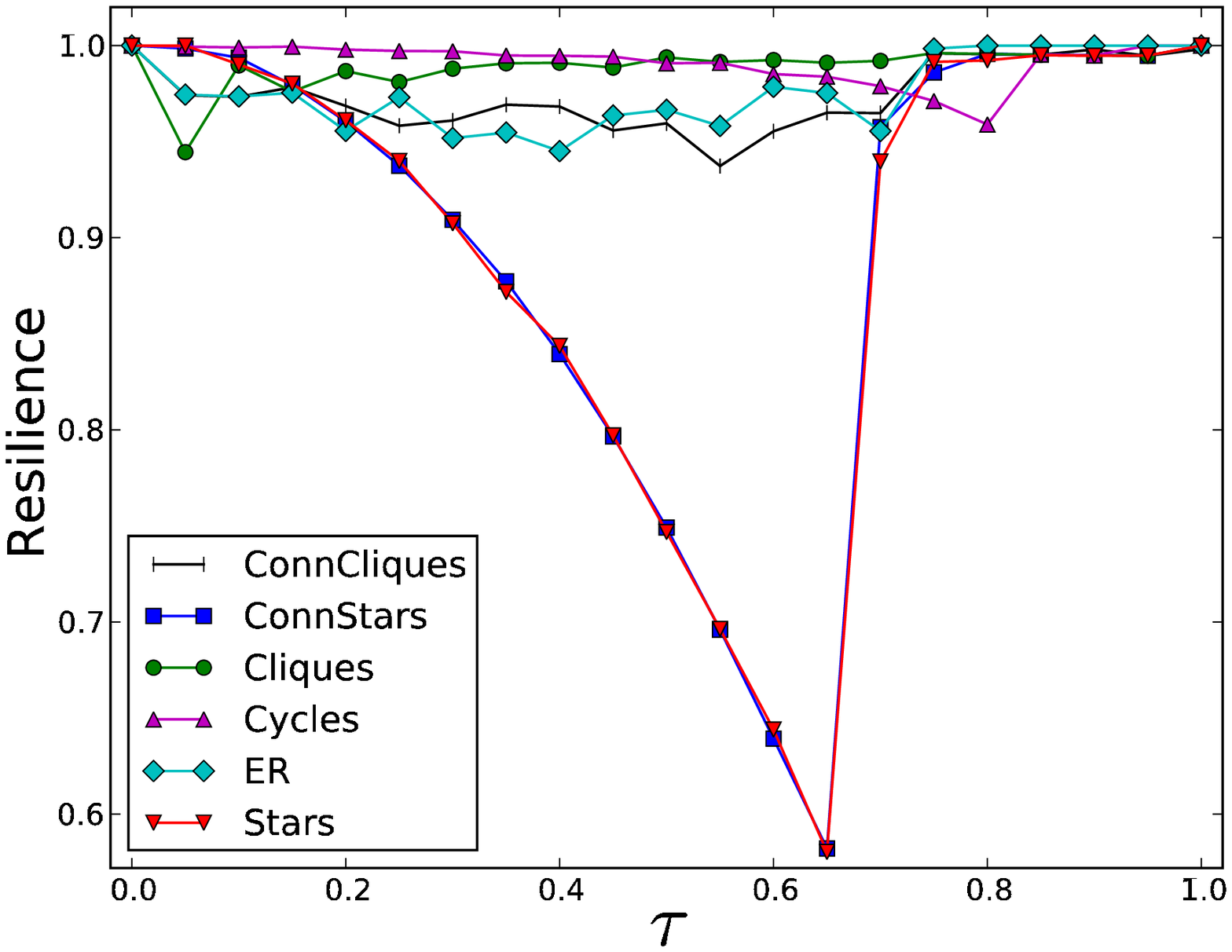}}
\par\end{centering}

\caption{{\bf Resilience of the optimal design.}}

\label{fig:opt-resilience-1}
\end{figure}

\begin{figure}[H]
\begin{centering}
\subfigure[$r=0.49$]{

\includegraphics[width=0.45\columnwidth]{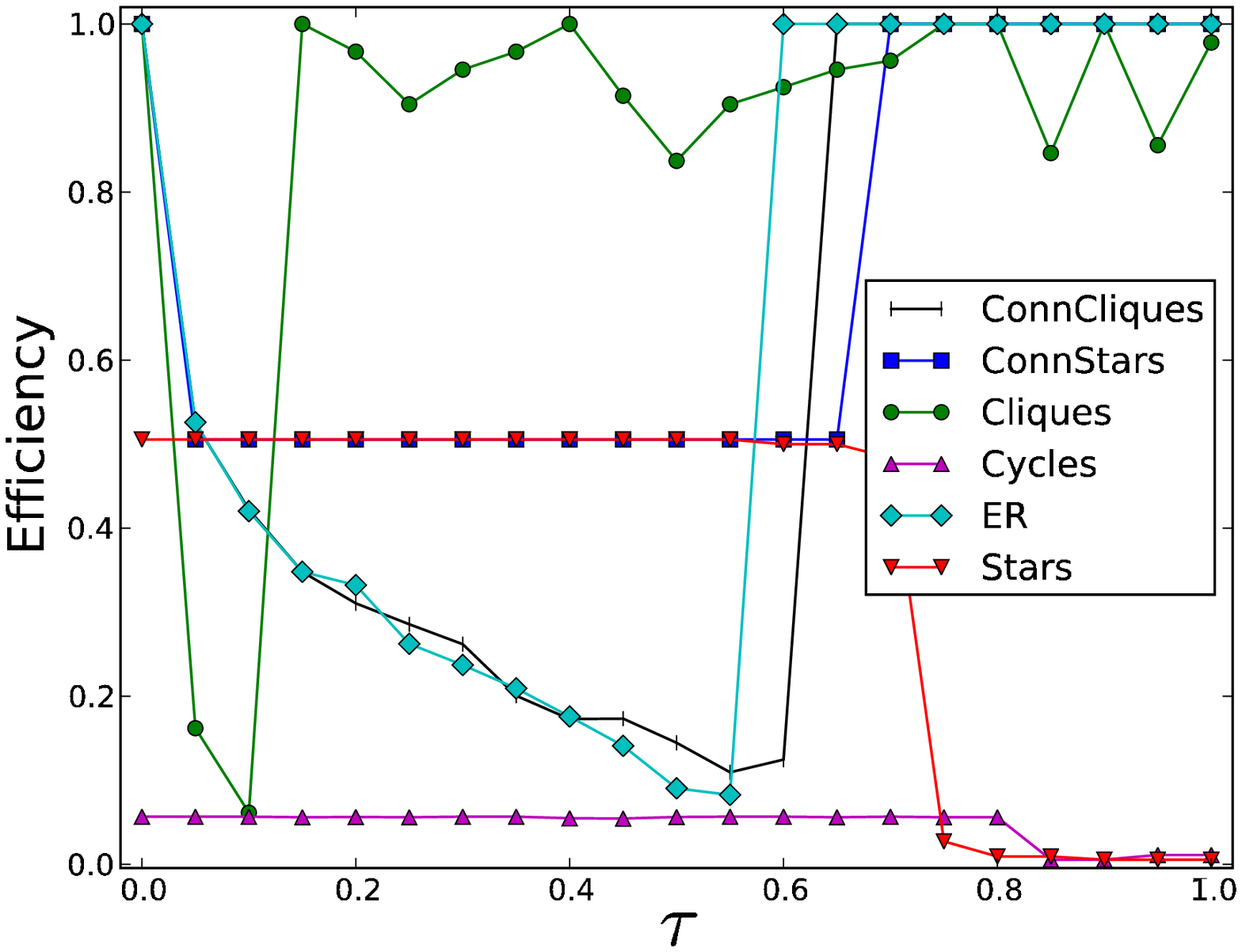}} \subfigure[$r=0.51$]{

\includegraphics[width=0.45\columnwidth]{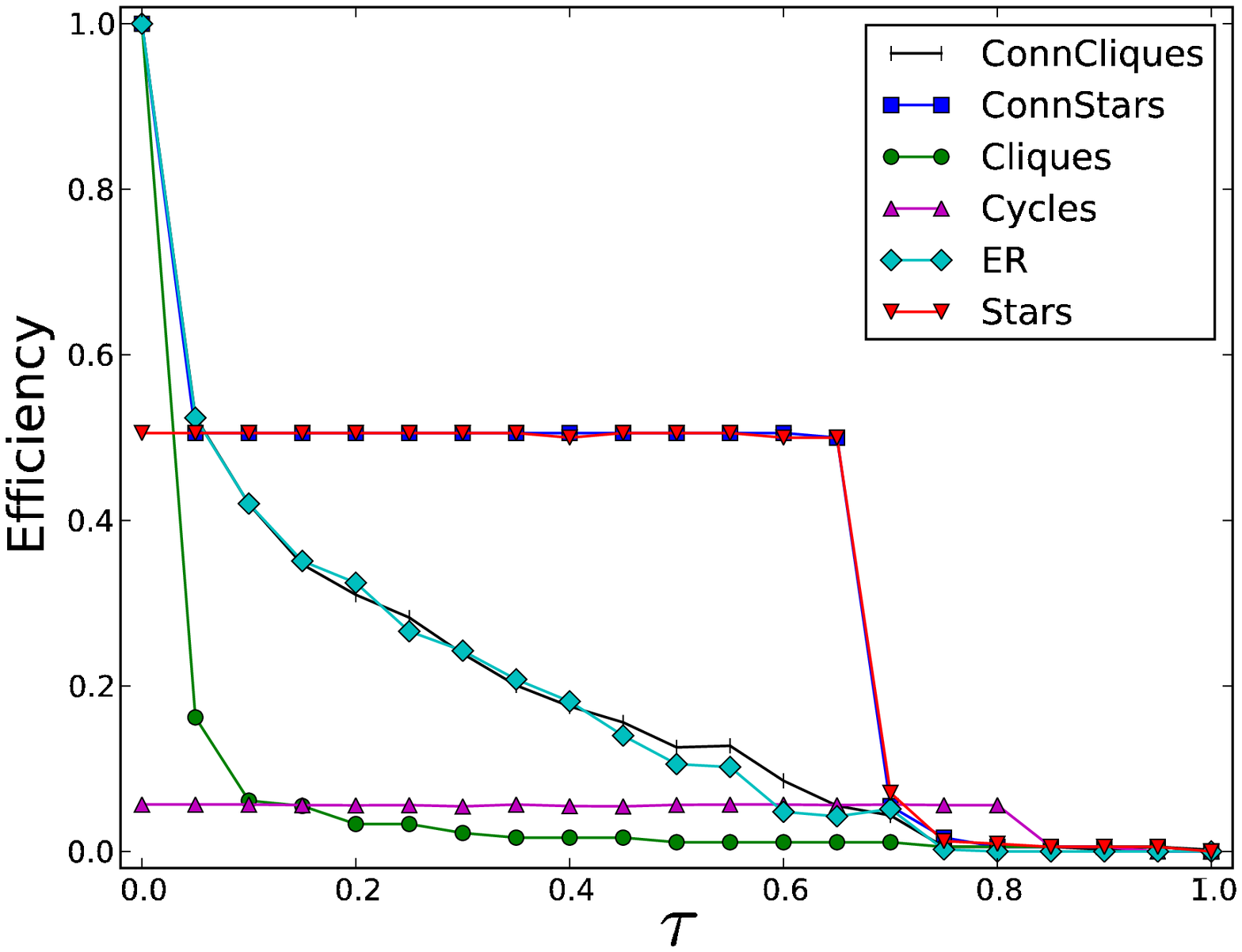}}
\par\end{centering}

\caption{{\bf Efficiency of the optimal design.}}

\label{fig:opt-efficiency-1}
\end{figure}

A plot of the fitness of various designs is in Fig.~\ref{fig:fitness-full}.
Notice that for $r=0.25$ Cliques outperform Stars. The resilience
and efficiency are in Fig.~\ref{fig:opt-resilience-1} and \ref{fig:opt-efficiency-1}. 

The success of the Stars design could be analyzed more qualitatively.
The fitness function combines resilience $R(G)$ which decreases when
the graph becomes more strongly connected, and efficiency $W(G)$
which decreases when the graph becomes more sparse. The optimality
of the star-based designs is due to a good trade-off between $R$
and $W$: the central node in each cell (its {}``leader'') provides
a good firewall against cascades because in each cell most pairs are
separated by distance of $2$, but this separation reduces efficiency
only modestly. In the Cliques design the separation is $1$ (too short
for resilience), and in the Cycles design it is too long ($\sim$a
quarter of cell size $k$).

Mathematically, the existence of a non-trivial solution is due to
the different functional relationships. To first-order approximation,
efficiency decreases inversely with average distance,$\Theta\left(\frac{1}{\mbox{avg distance}^{g}}\right)$
while cascade extent decreases exponentially, $\Theta\left(\tau^{\mbox{avg distance}}\right)$,
for $\tau<1$ assuming a bounded number of alternative paths). For
example, for the star design $R=O(1-\tau^{2})$ and $W=O(2^{-g})$
as $n=k\to\infty$. Therefore, the optimal network's structure exploits
the exponential decrease in cascades without sacrificing too much
efficiency. In the range $\tau\in[0.2,0.7]$ and $r\approx0.5$, an
average distance of $\approx2$, as in the star graph, might be optimal
(cf. \cite{Lindelauf08covert}).

\section{Monotonicity of Fitness \label{sec:Extent-and-Cascade}}

\emph{Proposition:} Let $f(\tau)$ be the highest attainable fitness
within a fixed network design $\mathbb{G}$, for cascade probability
$\tau$: \[
f(\tau)=\max_{G\in\mathbb{G}}\left[\underbrace{rR(G,\tau)+(1-r)W(G)}_{F(G,\tau)}\right]\]
Then $f(\tau)$ is a non-increasing function of $\tau$.

\emph{Proof of Proposition: }The proof relies on the simple claim
that resilience of networks does not increase when $\tau$ increases \cite{Gutfraind10extent}.
The claim is equivalent to the result that for a given graph $G$
increasing $\tau$ does not decrease the expected extent of cascades.
The remainder is almost trivial: we need to show that when the fitness
of all the points on the space (all graphs) has been made smaller
or kept the same (by increasing $\tau$), the new maximum value would
not be greater than in the old space. Rigorously, assume by contradiction
that $\tau_{+}>\tau$ and fitness \emph{increased}, namely\emph{:}
\begin{equation}
f(\tau_{+})>f(\tau)\,.\label{eq:fitness-contradict}\end{equation}
Let $G_{\tau_{+}}^{*},G_{\tau}^{*}$ by any two optimal networks for
$\tau_{+}$ and $\tau$, respectively, namely: 
\begin{align*}
G_{\tau_{+}}^{*}&\in\argmax_{G\in\mathbb{G}}\left[rR(G,\tau_{+})+(1-r)W(G)\right]\\
G_{\tau}^{*}&\in\argmax_{G\in\mathbb{G}}\left[rR(G,\tau)+(1-r)W(G)\right]\,.
\end{align*}
By optimality of $G_{\tau}^{*}$ at $\tau$ get that it must be at
least as good at $\tau$ as $G_{\tau_{+}}^{*}$: \begin{equation}
\underbrace{F(G_{\tau}^{*},\tau)-F(G_{\tau_{+}}^{*},\tau)}_{\equiv\Delta}\geq0\,.\label{eq:opt-tau}\end{equation}
The claim implies that:\begin{equation}
R(G_{\tau_{+}}^{*},\tau)\geq R(G_{\tau_{+}}^{*},\tau_{+})\,.\label{eq:claim}\end{equation}
Expanding $\Delta$:\begin{eqnarray*}
\Delta & = & F(G_{\tau}^{*},\tau)-\left[rR(G_{\tau_{+}}^{*},\tau)+(1-r)W(G_{\tau_{+}}^{*})\right]\\
 & \leq & F(G_{\tau}^{*},\tau)-rR(G_{\tau_{+}}^{*},\mathbf{\tau_{+}})-(1-r)W(G_{\tau_{+}}^{*})\mbox{ by }(\ref{eq:claim})\\
 & = & F(G_{\tau}^{*},\tau)-F(G_{\tau_{+}}^{*},\tau_{+})\\
 & < & 0\mbox{ by the hypothesis (\ref{eq:fitness-contradict}).}\end{eqnarray*}
This contradicts Ineq.\ \ref{eq:opt-tau}. The argument is easy to
generalize. One could apply this method to the parameter $g$ of attenuation,
showing that fitness is non-increasing when attenuation is increased.

\section{Analytic Results\label{sec:Analytic-Results}}

It is easy to analytically derive the values of the resilience, efficiency
(and hence fitness) functions for certain simple designs: the Cycles
and the Stars designs. These are useful to gain deeper insight into
the effect of parameters. Recall that $n$ is the number of nodes
and $k$ is the number of nodes per cell. For $k=1$, in both designs
$R=1$ and $W=0$. When $k\geq2$, it is easy to verify that for the
Cycle design:\begin{eqnarray*}
R(n,k,\tau) & = & 1-\frac{1}{n-1}\left[2\tau\frac{1-\tau^{k-1}}{1-\tau}-(k-1)\tau^{k}\right]\\
W(n,k,g) & = & \frac{1}{n-1}\begin{cases}
\frac{1}{(\frac{k}{2})^{g}}+2\sum_{j=1}^{\frac{k}{2}-1}\frac{1}{j^{g}} & \mbox{ \mbox{k }even}\\
2\sum_{j=1}^{\frac{k-1}{2}}\frac{1}{j^{g}} & \mbox{ \mbox{k }odd}\end{cases}\end{eqnarray*}
and for the Stars design:

\begin{eqnarray*}
R(n,k,\tau) & = & 1-\frac{1-\frac{1}{k}}{n-1}\left[2+\tau(k-2)\right]\tau\\
W(n,k,g) & = & \frac{1-\frac{1}{k}}{n-1}\left[2+2^{-g}(k-2)\right]\,.\end{eqnarray*}
These expressions are not readily useful for continuous optimization
since $k$ is discrete, but they can be used to identify phase transitions.
Thus, they help inform optimization for designs where no analytic
expression is available. The ER and Cliques designs are also analytically
tractable (see e.g. \cite{Draief08} for expected cascade extent).

\begin{figure}[H]
\begin{centering}
\subfigure[$r=0.25$]{

\includegraphics[width=0.35\columnwidth]{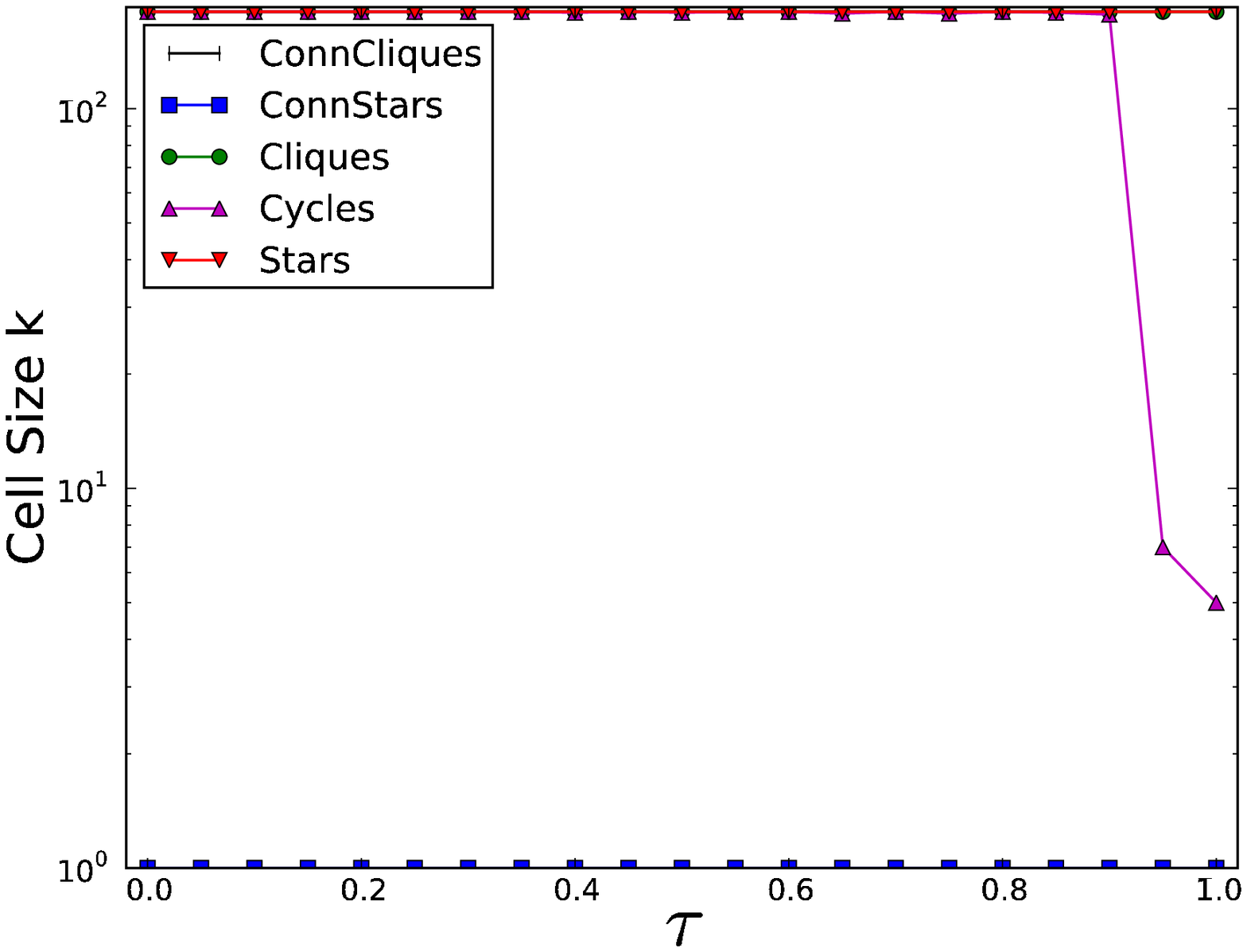}} \subfigure[$r=0.49$]{

\includegraphics[width=0.35\columnwidth]{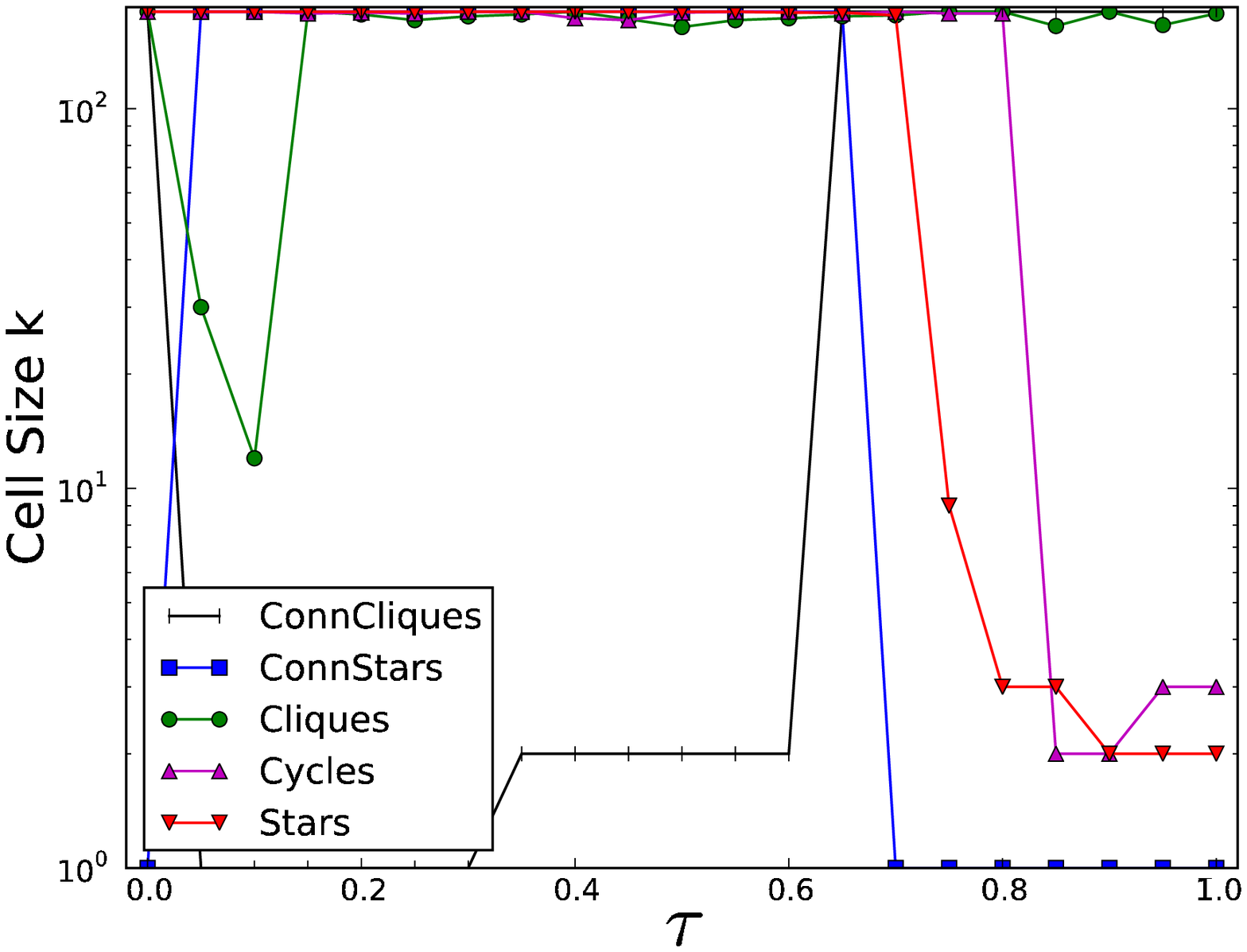}}
\par\end{centering}

\begin{centering}
\subfigure[$r=0.51$]{

\includegraphics[width=0.35\columnwidth]{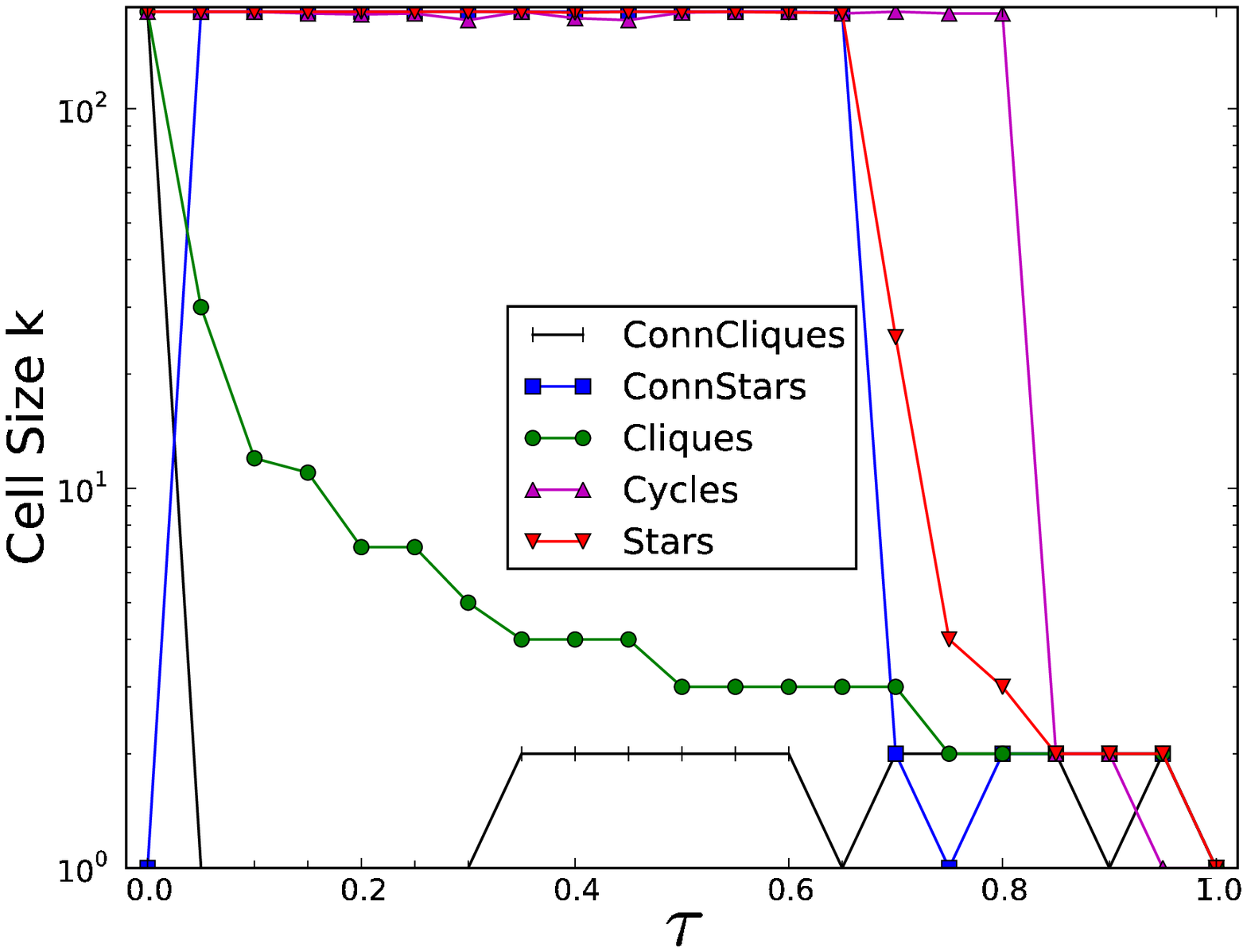}} \subfigure[$r=0.75$]{

\includegraphics[width=0.35\columnwidth]{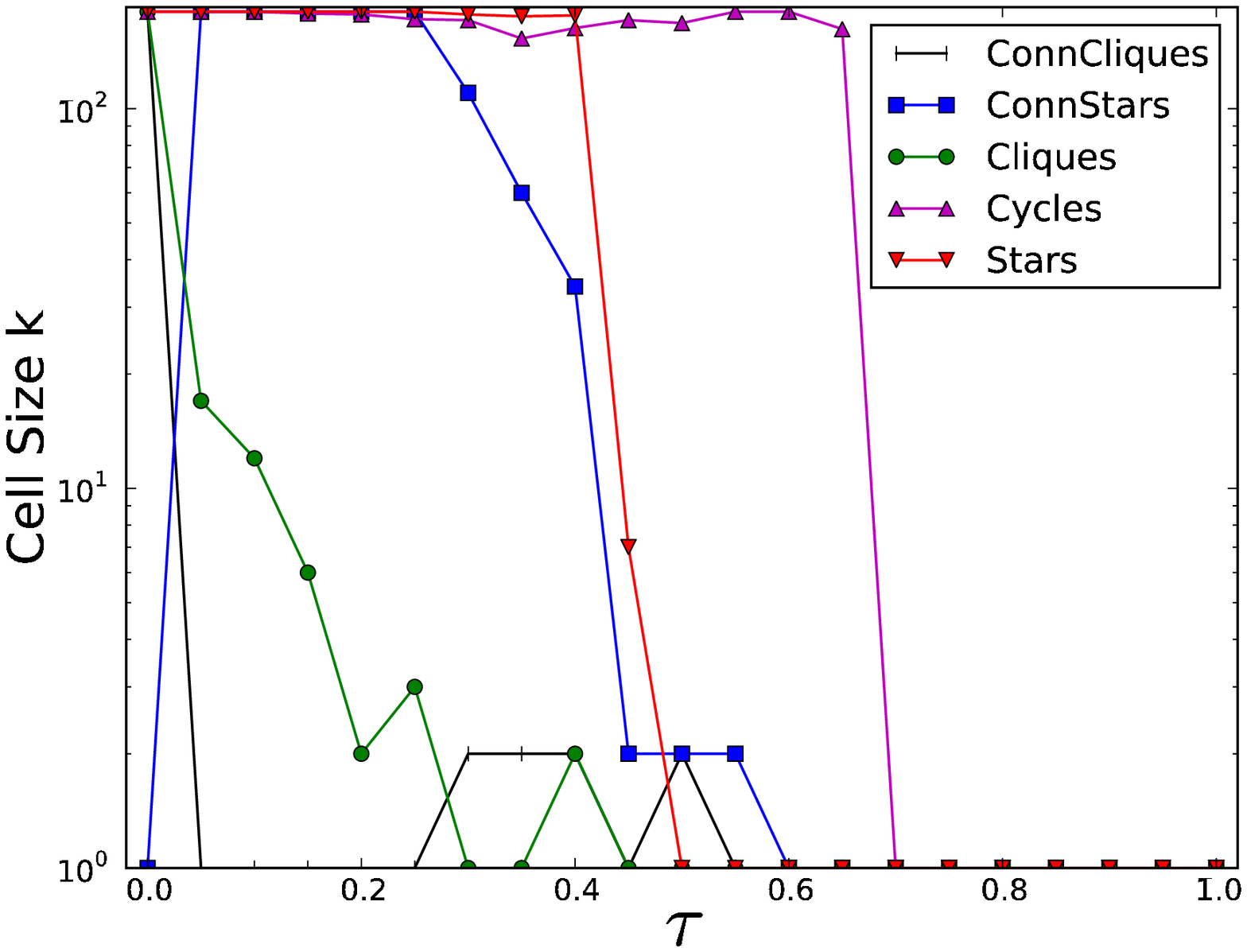}}
\par\end{centering}

\caption{{\bf Cell size $k$ in the optimal configuration of each design.} Cell size makes many sharp transitions.}

\label{fig:opt-k}
\end{figure}

In the Stars design, when $R$ and $W$ are weighted equally ($r=\frac{1}{2}$),
fitness takes a relatively simple form: $F=\frac{1}{2}+\frac{1}{2}\frac{1-\frac{1}{k}}{n-1}\left[2(1-\tau)+(2^{-g}-\tau^{2})(k-2)\right]$.
This implies that increasing cell size $k$, for $k$ large, improves
fitness iff $2^{-g}-\tau^{2}>0$. Hence the optimal configuration
has one cell ($k=n$), until a threshold near $\tau=2^{-g/2}$ (for
$g=1$, approximately $0.71$). This agrees with the findings in Fig.~\ref{fig:opt-k}.
Also, the rate of change in fitness with respect to $\tau$, $\frac{dF}{d\tau}=\frac{1-\frac{1}{k}}{n-1}\left[-2-2\tau(k-2)\right]$,
is always negative, as expected on more general grounds (see sec.\,\ref{sec:Extent-and-Cascade}). It is linear in $\tau$ (because
it is a tree graph) but superlinear in $k$ (because of the mutual
hazard induced by adding nodes to cells.)

\section{Configuration of the Optimal Design\label{sec:Configuring-the-Optimal}}

As $\tau$ is varied, the optimal configuration changes. This section
shows the resulting changes in the values of the parameters $k$ (cell
size) and $p$ (connectivity). In other words, it indicates how each
of the designs ought to be configured to attain optimal fitness, as
a function of resilience weighting, $r$, and cascade probability,
$\tau$.

The cell size parameter $k$ is non-monotonic for various designs
when $r=0.49$ (Fig.~\ref{fig:opt-k}). For example, for the Connected
Cliques design, at low contagion risk ($\tau<0.1$), $k$ is high
(comparable to the size of the network, i.e. $k\to n$), then it falls
to a small number. At high contagion risk ($\tau>0.6$) the network
is again highly connected again with $k\to n$. Thus for $\tau\to1$,
the optimal network is the fully-connected graph. 

\begin{figure}[H]
\begin{centering}
\subfigure[$r=0.25$]{

\includegraphics[width=0.35\columnwidth]{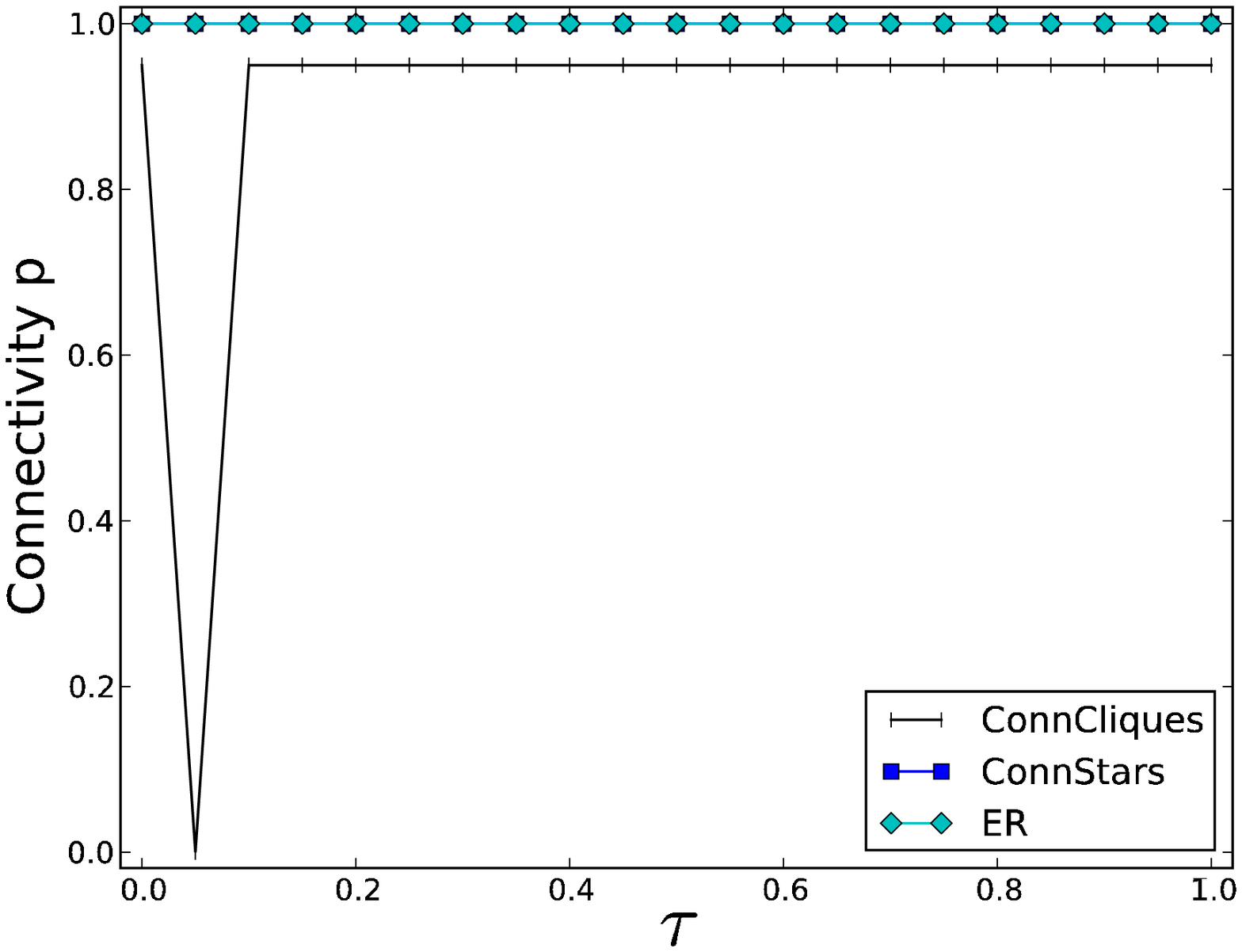}} \subfigure[$r=0.49$]{

\includegraphics[width=0.35\columnwidth]{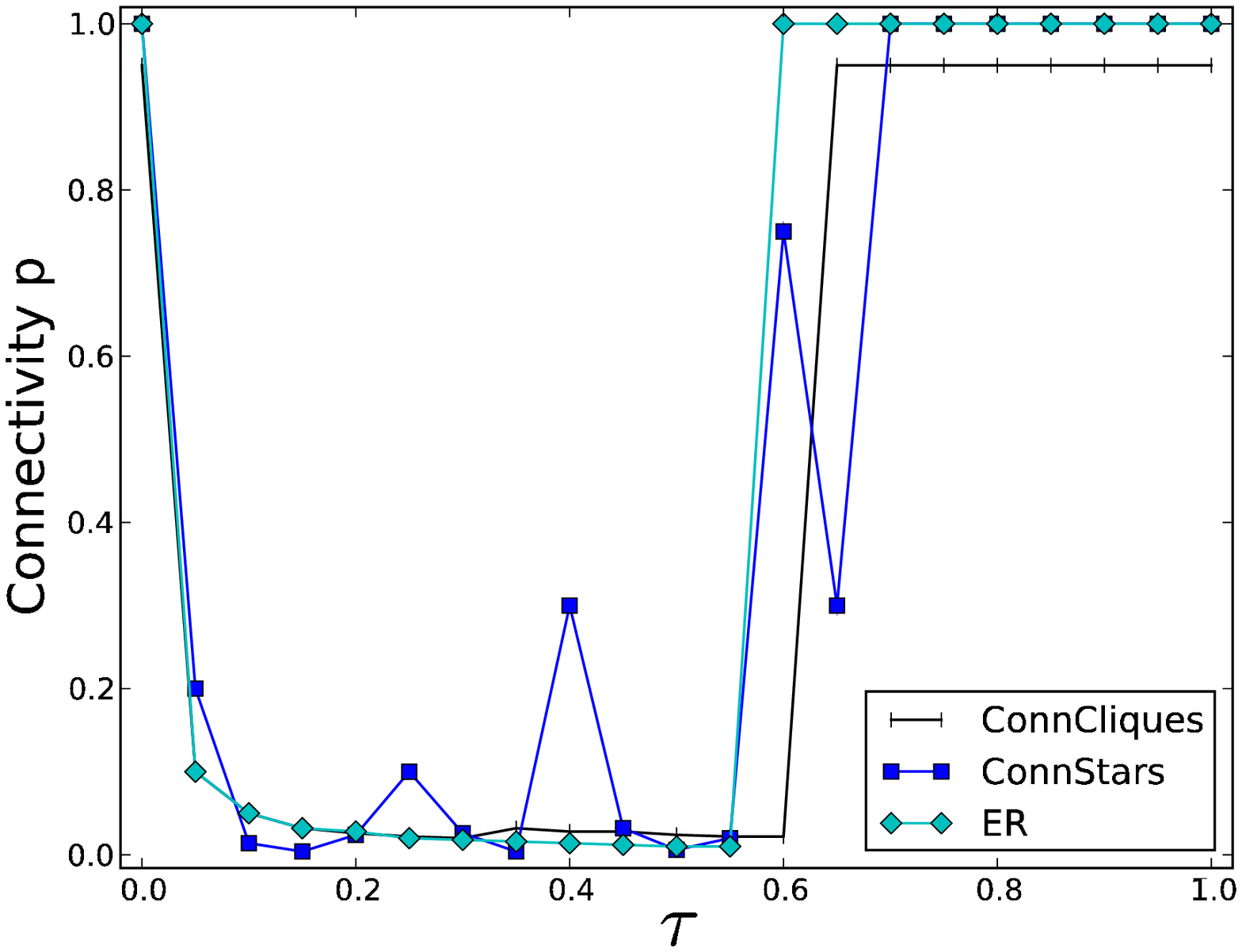}}
\par\end{centering}

\begin{centering}
\subfigure[$r=0.51$]{

\includegraphics[width=0.35\columnwidth]{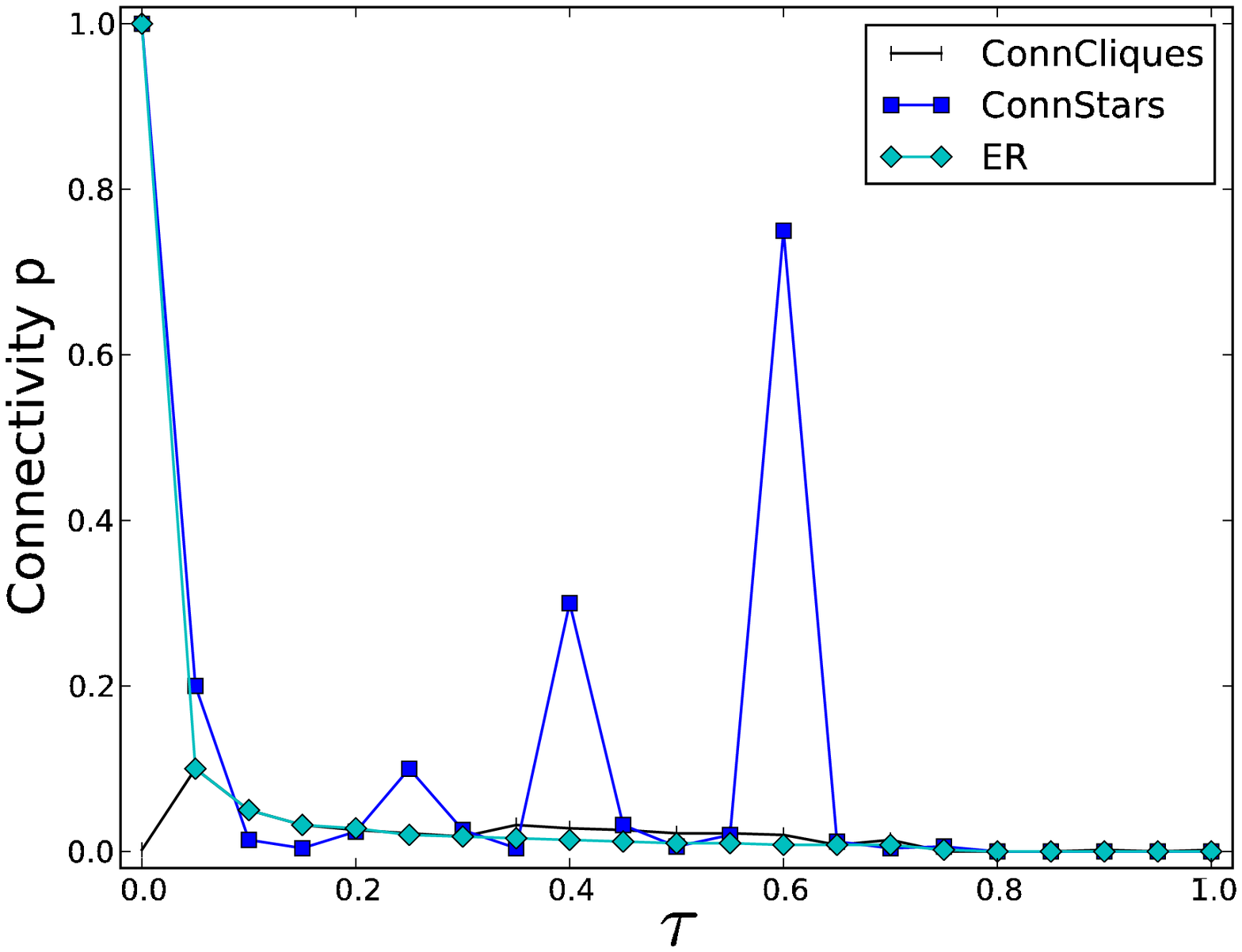}} \subfigure[$r=0.75$]{

\includegraphics[width=0.35\columnwidth]{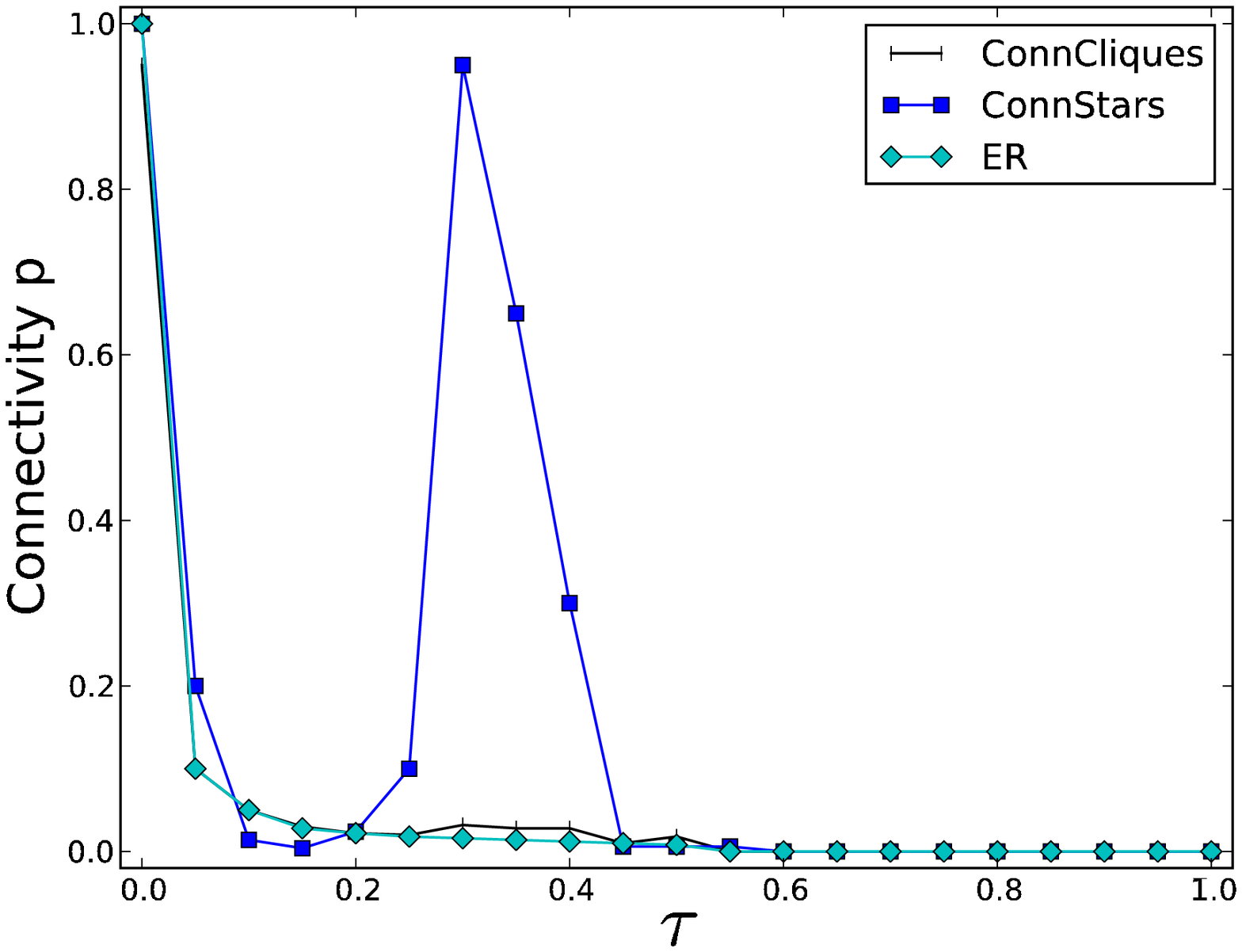}}
\par\end{centering}

\caption{{\bf Connectivity $p$ in the optimal configuration of each design.} Non-monotonic trends reflect structural transitions and sometimes also interplay with cell size.}

\label{fig:opt-p}
\end{figure}

In general, designs involving both the $p$ and $k$ parameters show
an interplay between the two (Fig.~\ref{fig:opt-p}). For example,
in the Connected Stars design under $r<0.5$ there are two phase-transitions
in connectivity $p$: as $\tau$ increases at $\tau\to\tau_{l}^{*}\approx0.1$
it transitions from a connected graph to disconnected cells, and at
$\tau\to\tau_{u}^{*}\approx0.7$ back to full connectivity. If $r>0.5$
the second transition is extinguished. The data requires care to interpret.
For example, in the Connected Stars design $\tau\in[0.1,0.65]$, when
$r>0.5$ the fluctuations in the $p$ are noise because there is a
single cell and a single cell leader ($k=n$), and so the parameter
$p$ has no effect. For sensitivity analysis see section~\ref{sec:Sensitivity-Analysis}.\\

\section{Sensitivity Analysis\label{sec:Sensitivity-Analysis}}

It is desirable to determine how much variability exists within the
optimal parameter values. Only if the variability is low the optimal
configuration provides valuable information for network design. As
proxy of robustness of the estimates, consider the space of configurations
whose fitness $\geq0.95$ of the fitness of the optimal solution.
The robustness of a parameter is measured by the parameter's standard
deviation within this space (the unbiased estimate). 

\begin{figure}[H]
\begin{centering}
\subfigure[$r=0.25$]{

\includegraphics[width=0.30\columnwidth]{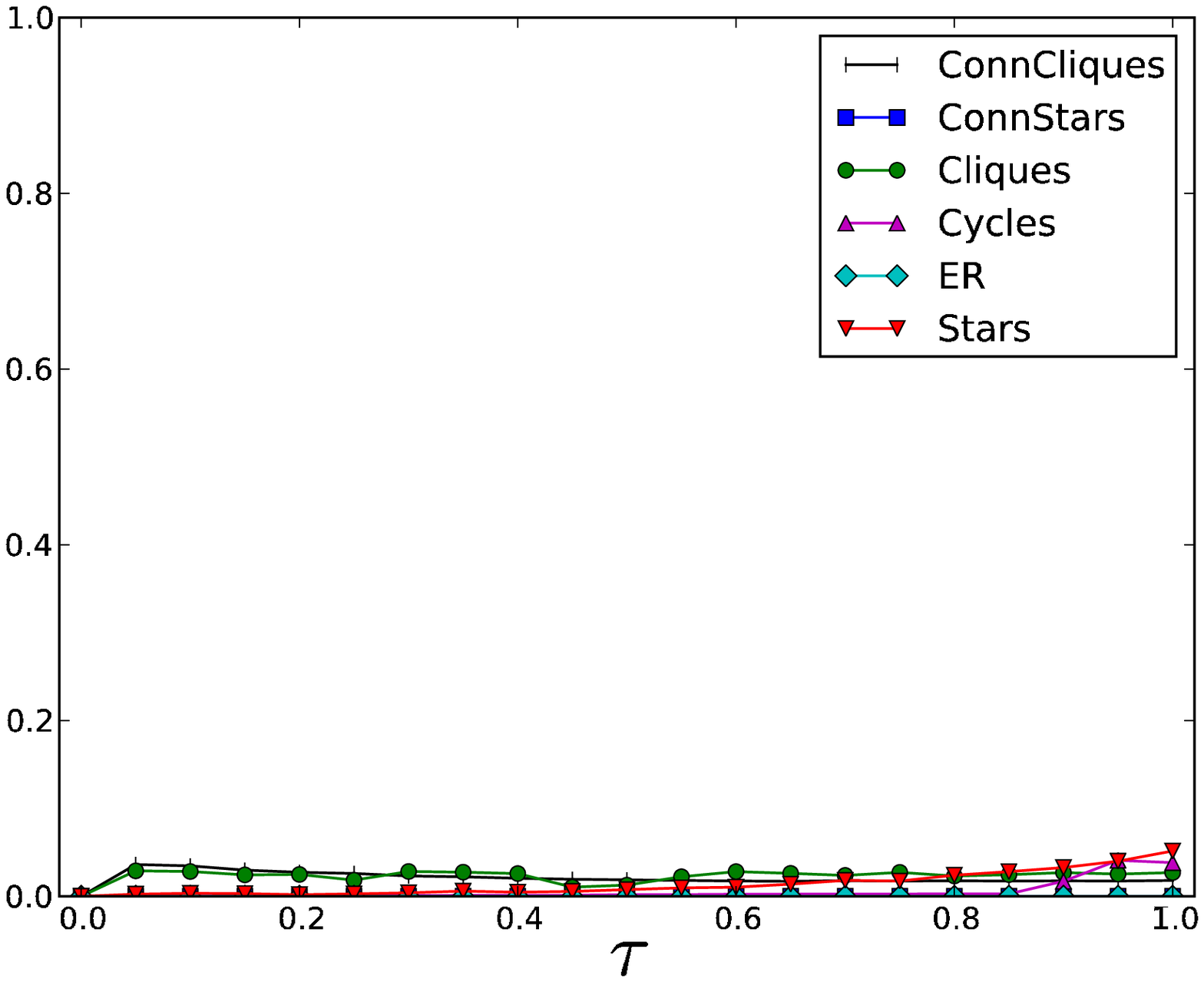}} \subfigure[$r=0.49$]{

\includegraphics[width=0.30\columnwidth]{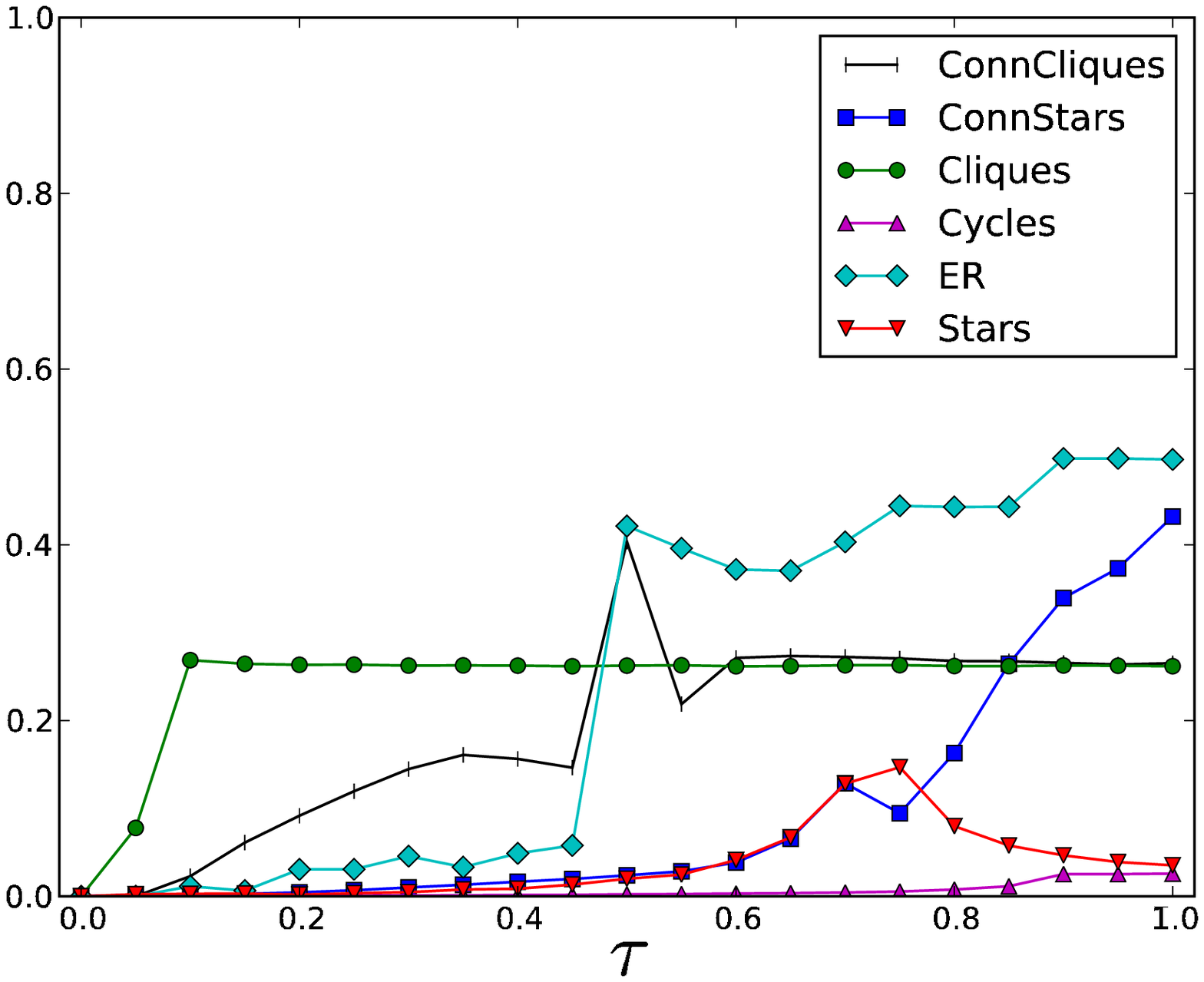}}
\par\end{centering}

\begin{centering}
\subfigure[$r=0.51$]{

\includegraphics[width=0.30\columnwidth]{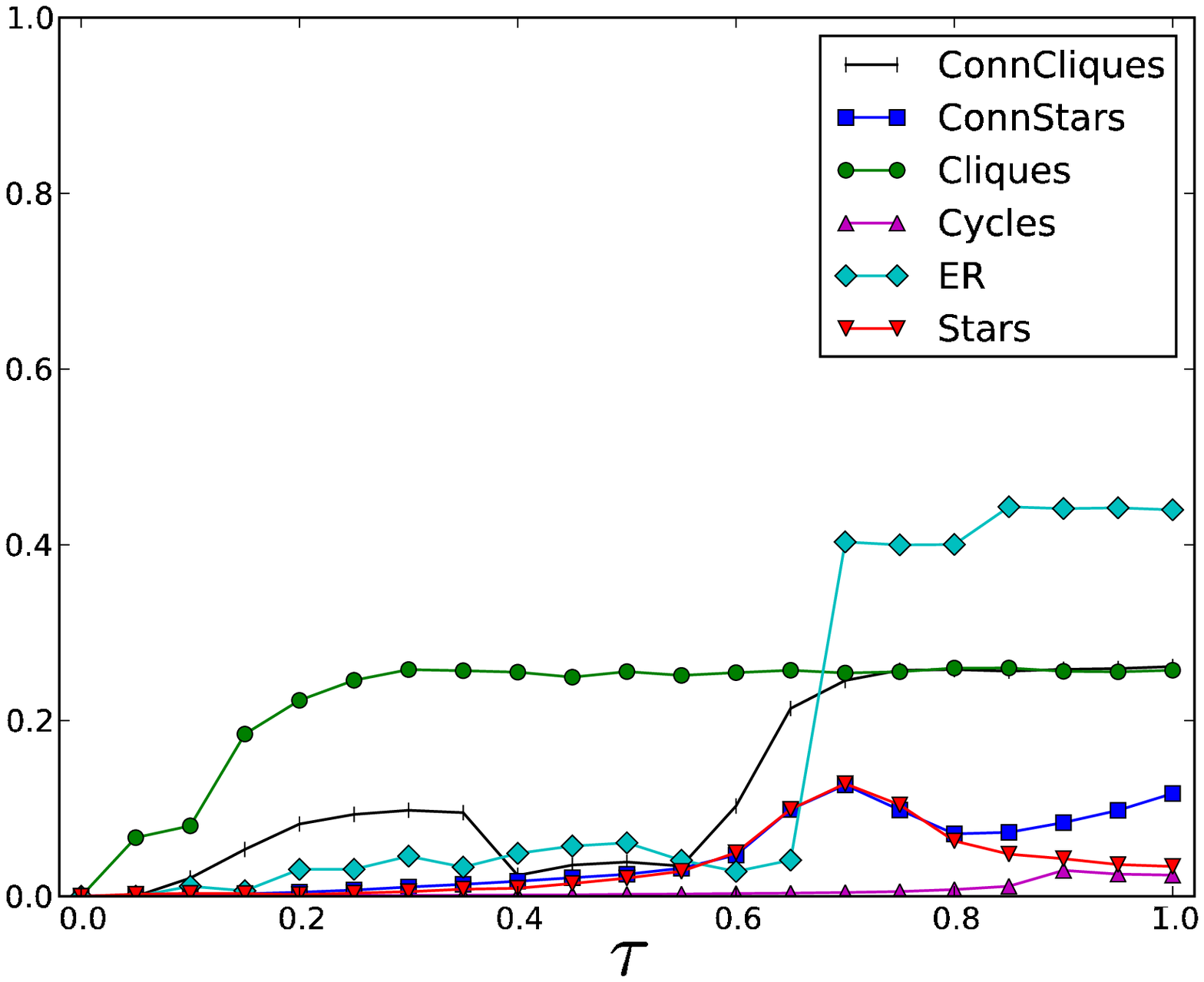}} \subfigure[$r=0.75$]{

\includegraphics[width=0.30\columnwidth]{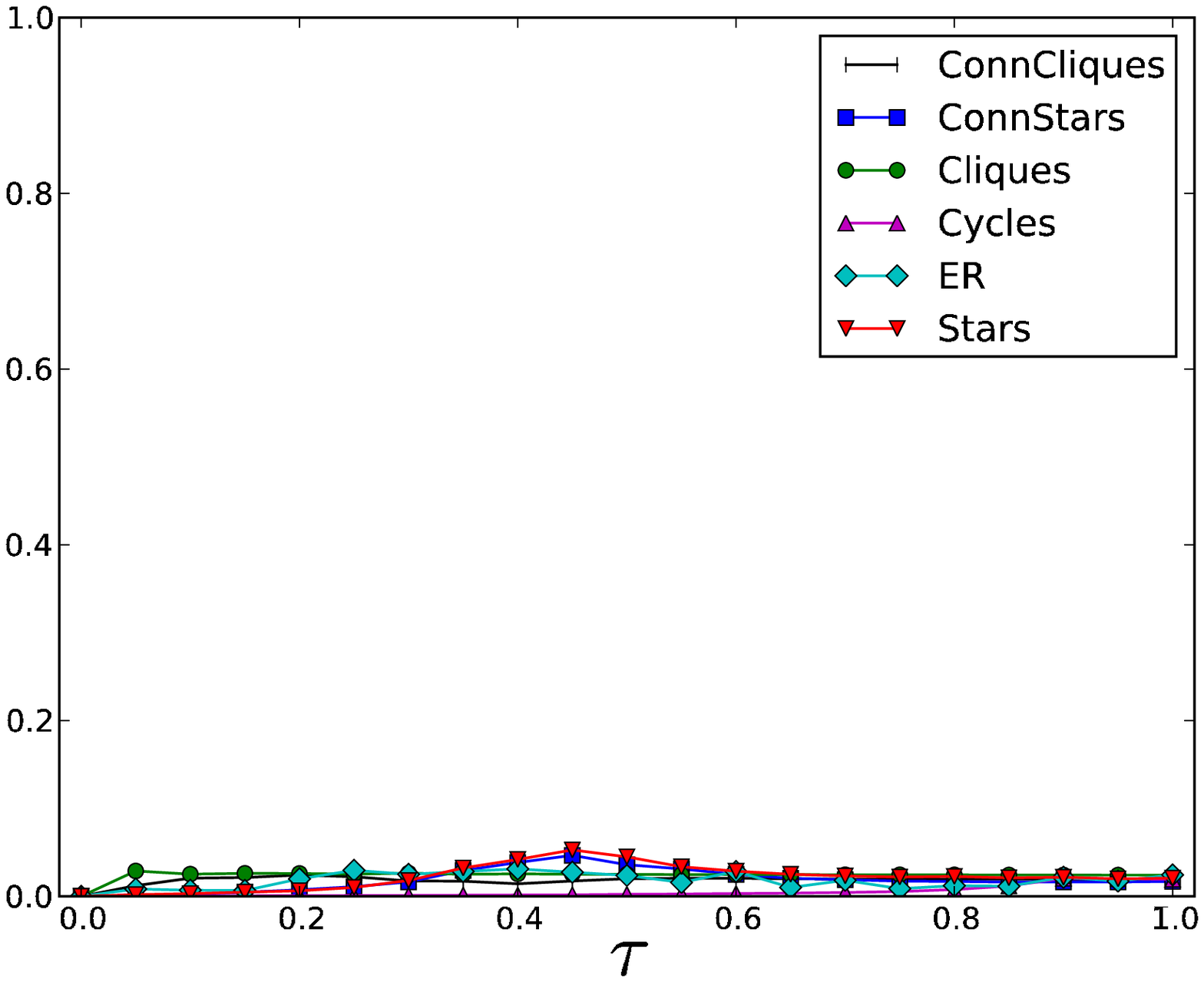}}
\par\end{centering}

\caption{{\bf Standard deviation in resilience within the top 5\% of solutions.}}
\label{fig:robustness-opt-resilience}
\end{figure}

\begin{figure}[H]
\begin{centering}
\subfigure[$r=0.25$]{

\includegraphics[width=0.30\columnwidth]{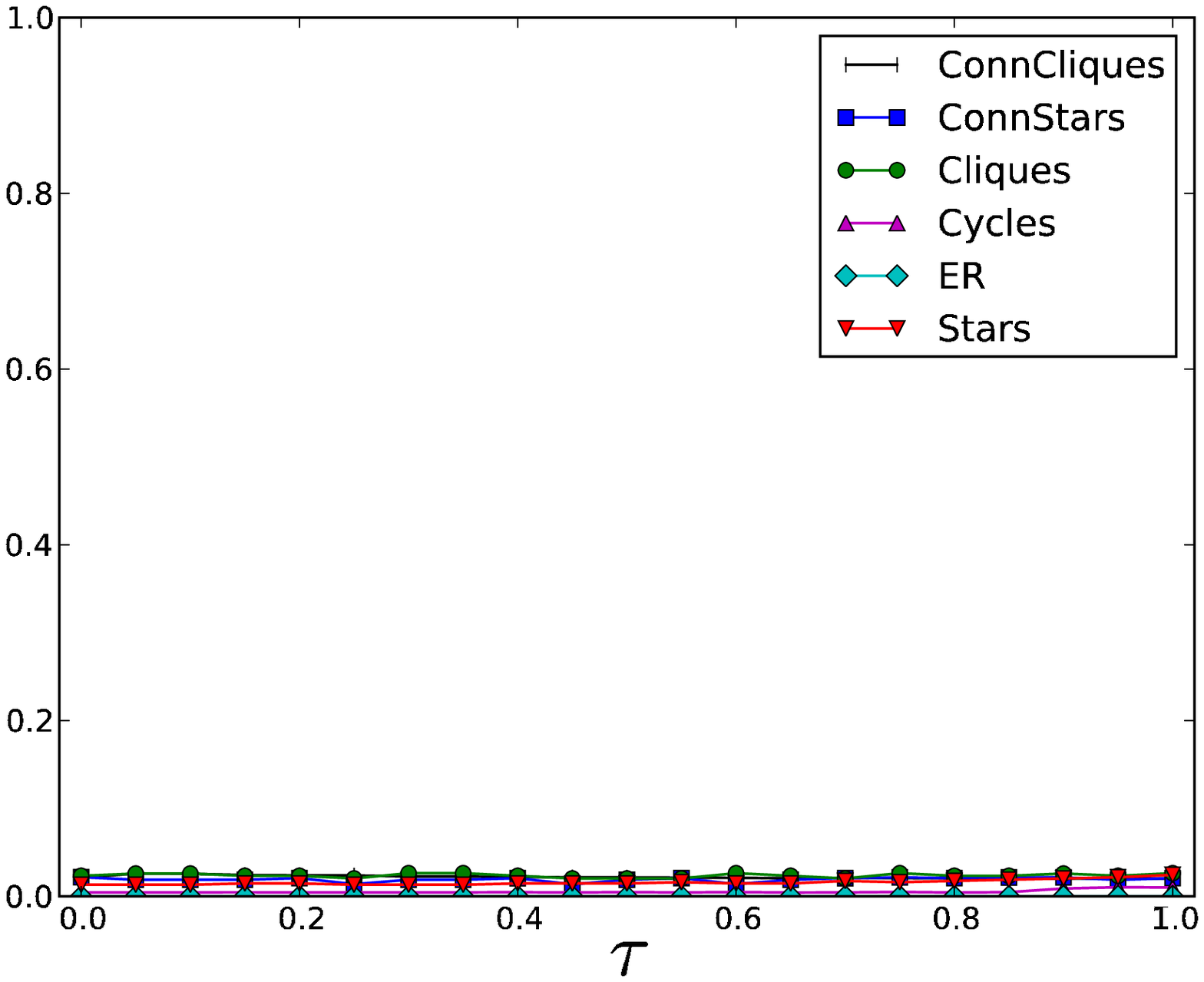}} \subfigure[$r=0.49$]{

\includegraphics[width=0.30\columnwidth]{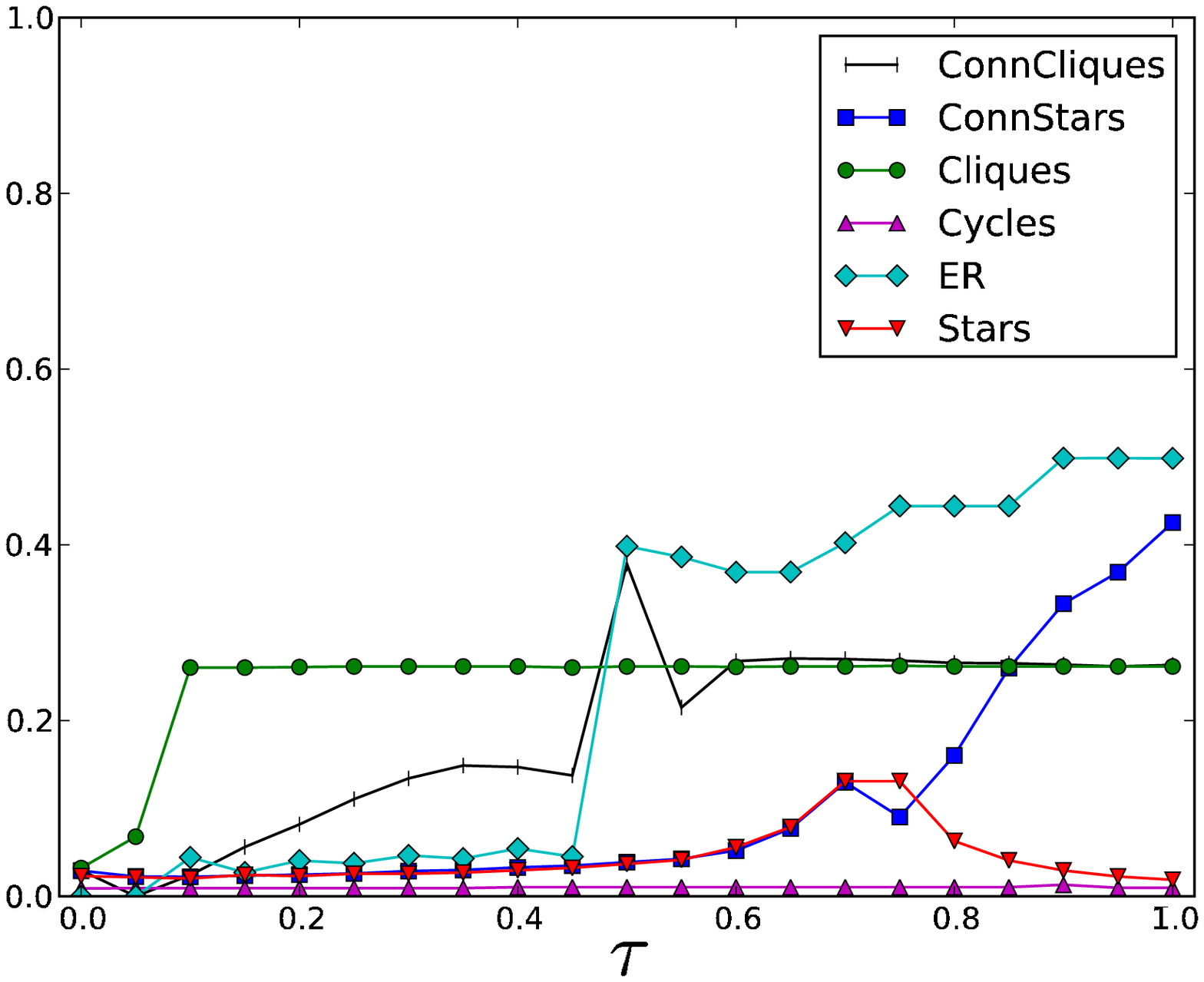}}
\par\end{centering}

\begin{centering}
\subfigure[$r=0.51$]{

\includegraphics[width=0.30\columnwidth]{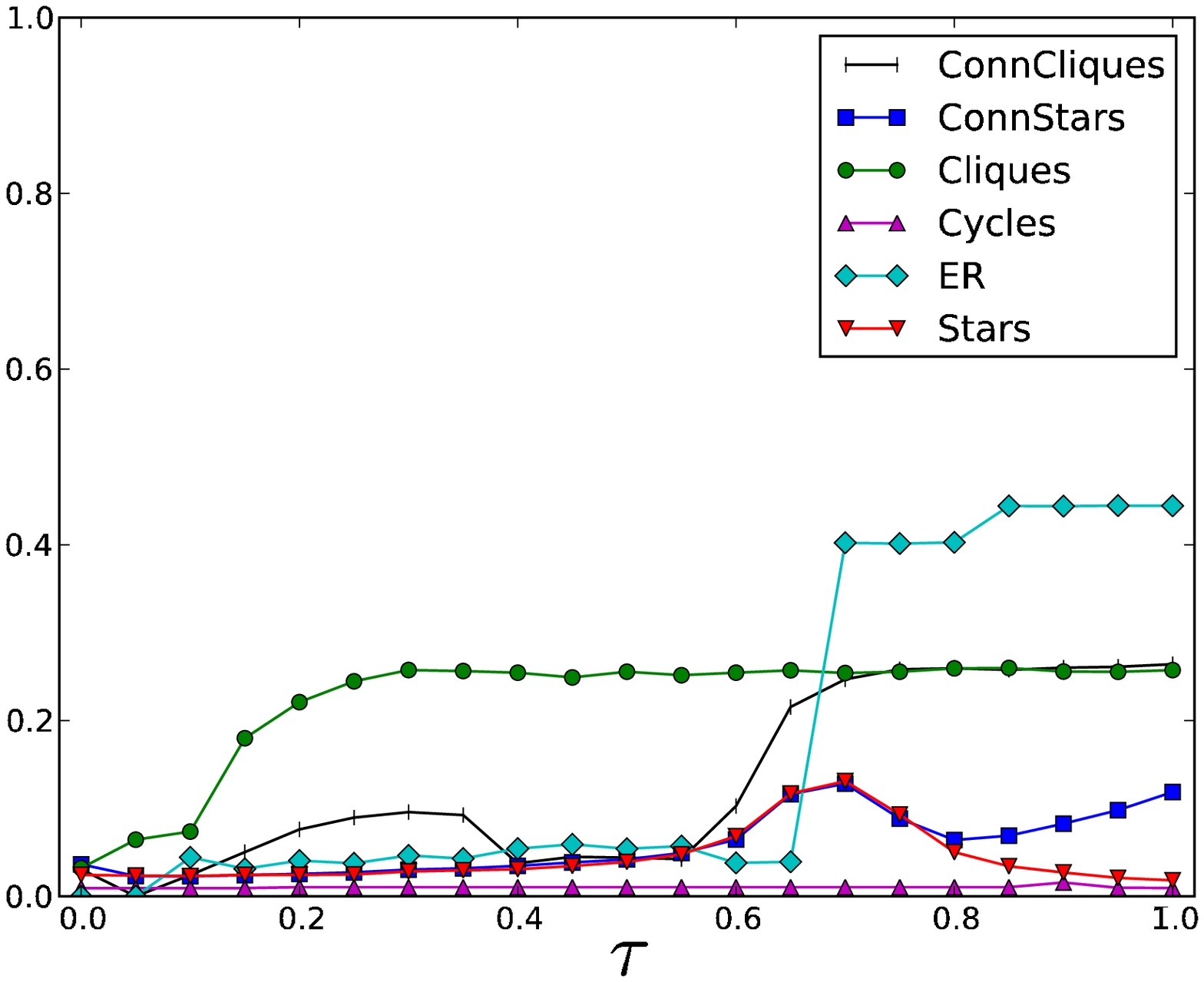}} \subfigure[$r=0.75$]{

\includegraphics[width=0.30\columnwidth]{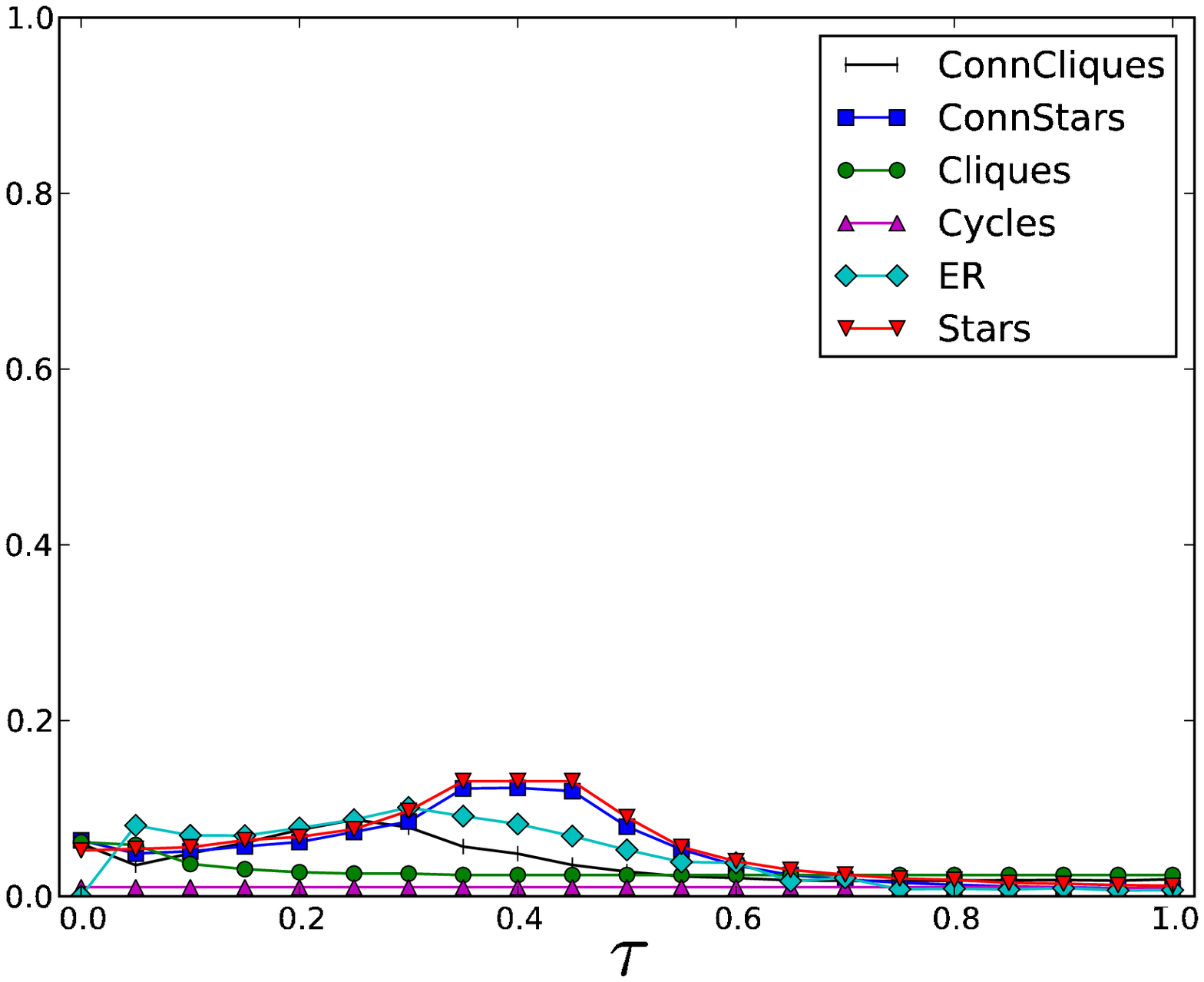}}
\par\end{centering}

\caption{{\bf Standard deviation in efficiency within the top 5\% of solutions.}}
\label{fig:robustness-opt-efficiency}
\end{figure}

\begin{figure}[H]
\begin{centering}
\subfigure[$r=0.25$]{

\includegraphics[width=0.30\columnwidth]{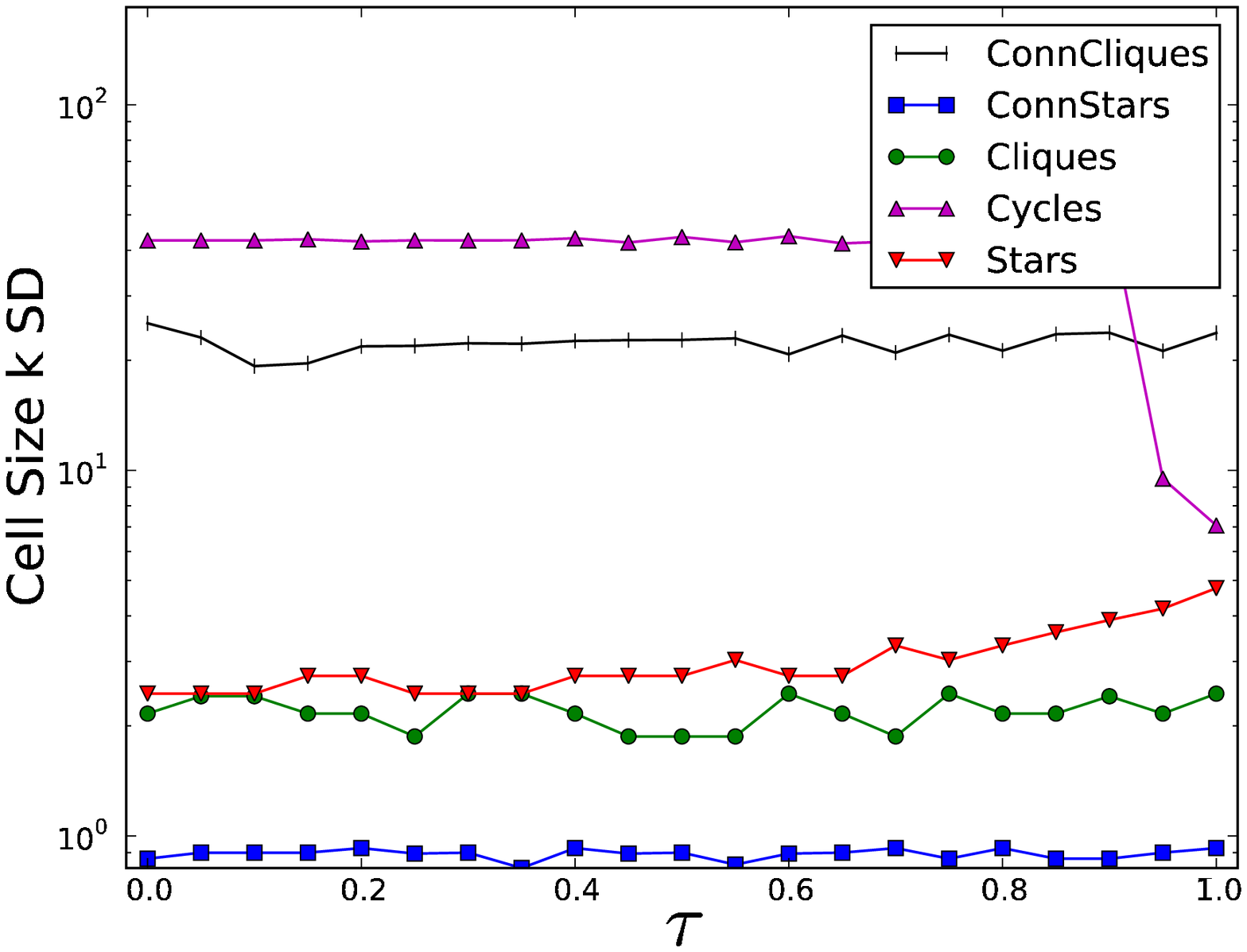}} \subfigure[$r=0.49$]{

\includegraphics[width=0.30\columnwidth]{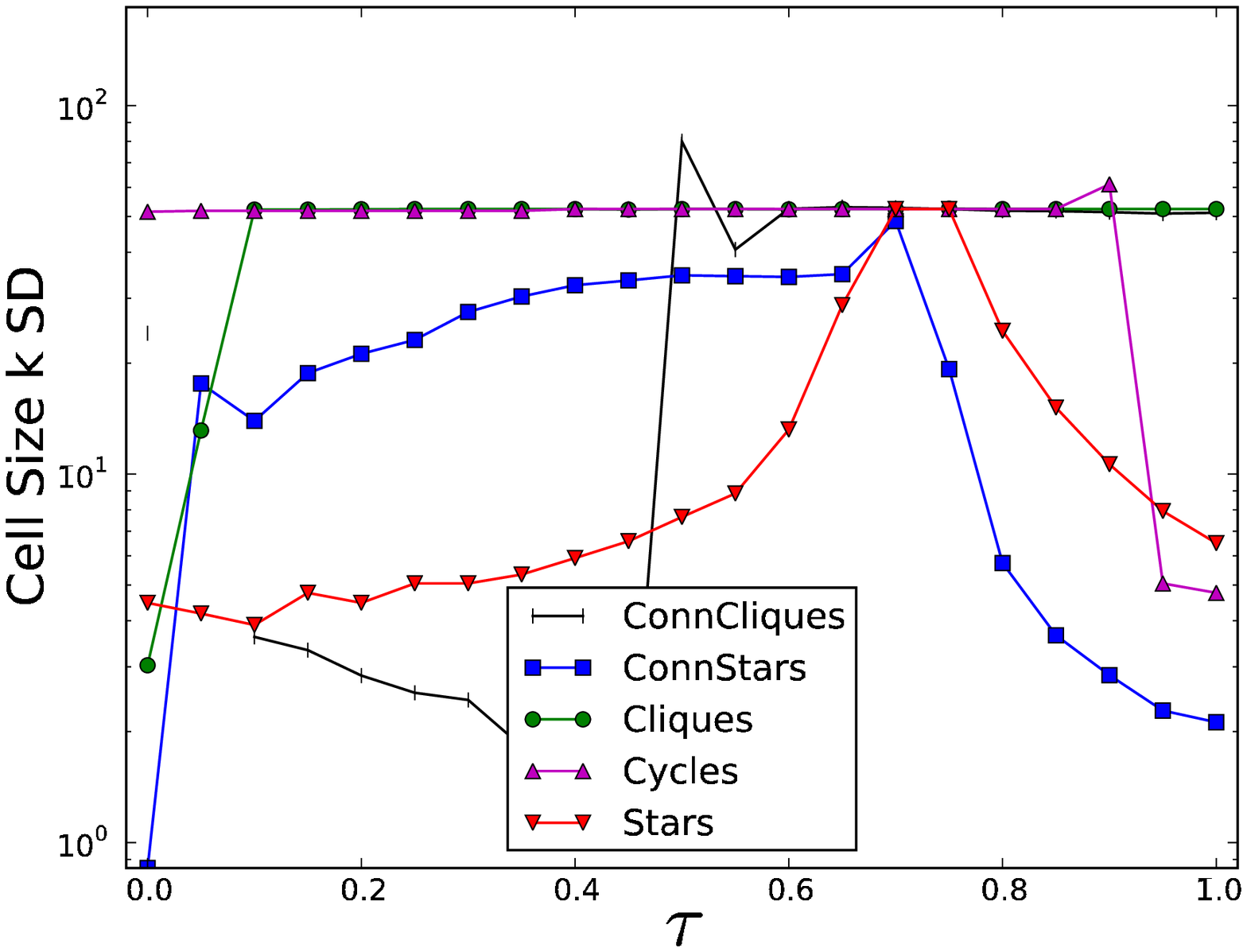}}
\par\end{centering}

\begin{centering}
\subfigure[$r=0.51$]{

\includegraphics[width=0.30\columnwidth]{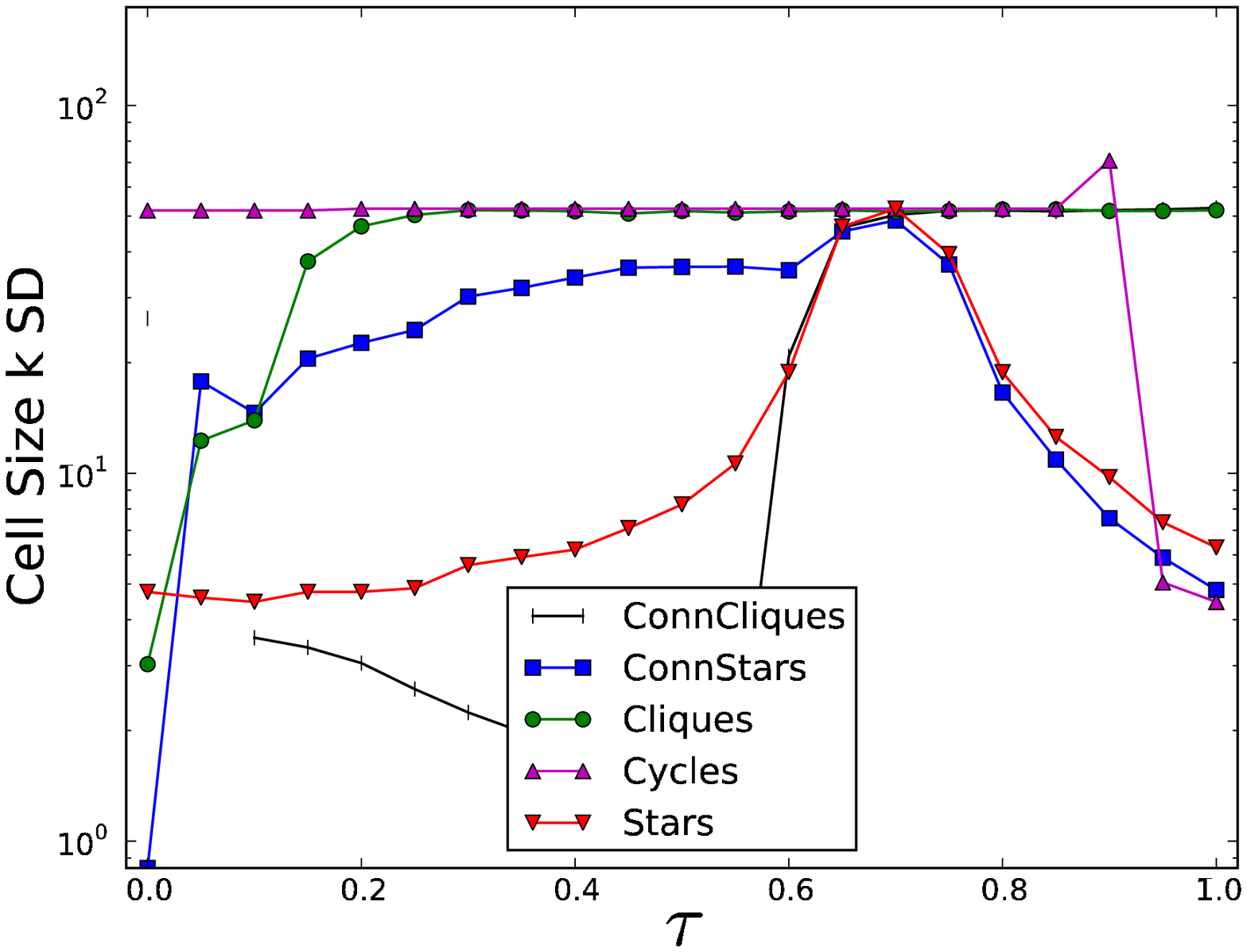}} \subfigure[$r=0.75$]{

\includegraphics[width=0.30\columnwidth]{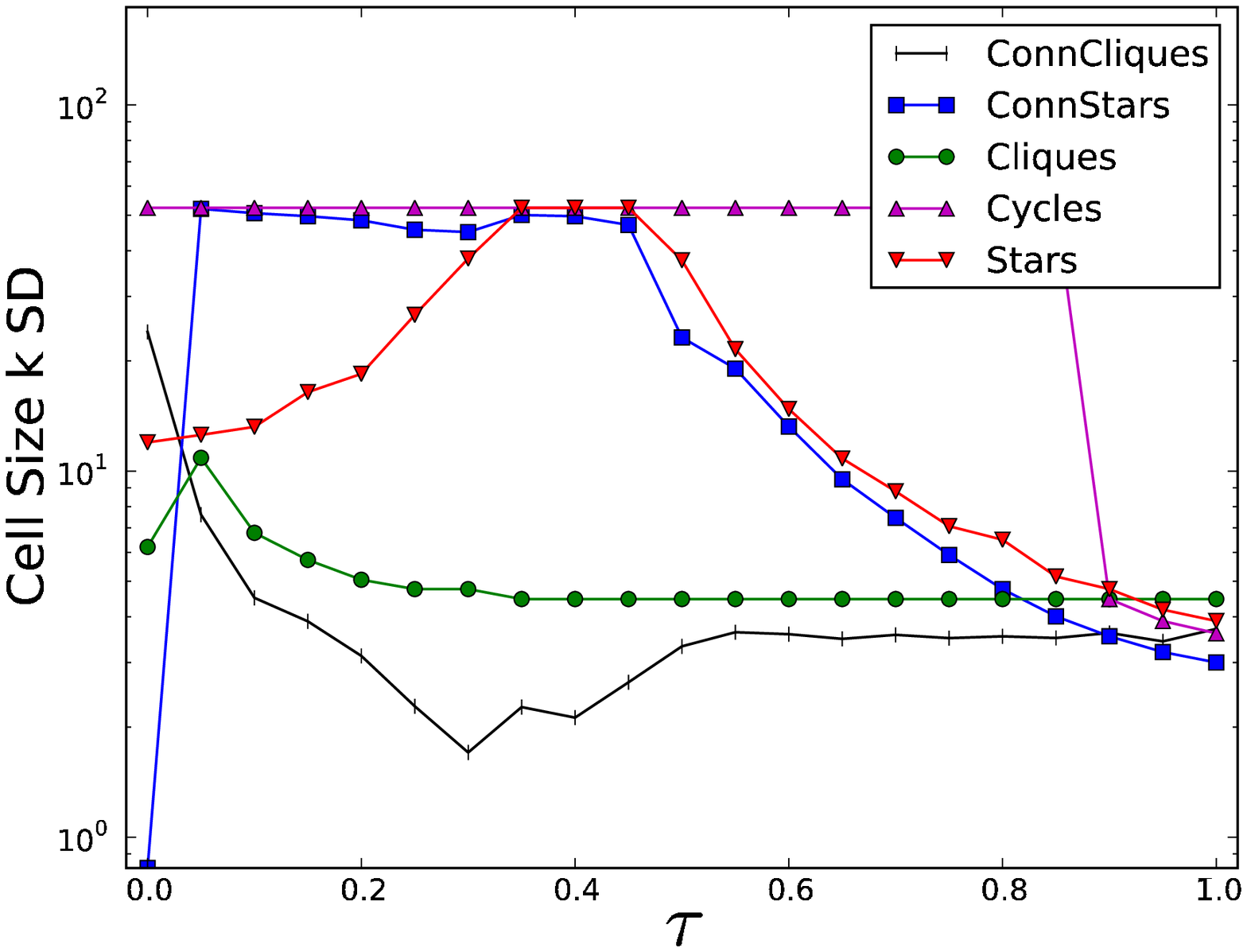}}
\par\end{centering}

\caption{{\bf Standard deviation in cell size $k$ of the top 5\% of solutions.}}
\label{fig:robustness-opt-k}
\end{figure}

\begin{figure}[H]
\begin{centering}
\subfigure[$r=0.25$]{

\includegraphics[width=0.30\columnwidth]{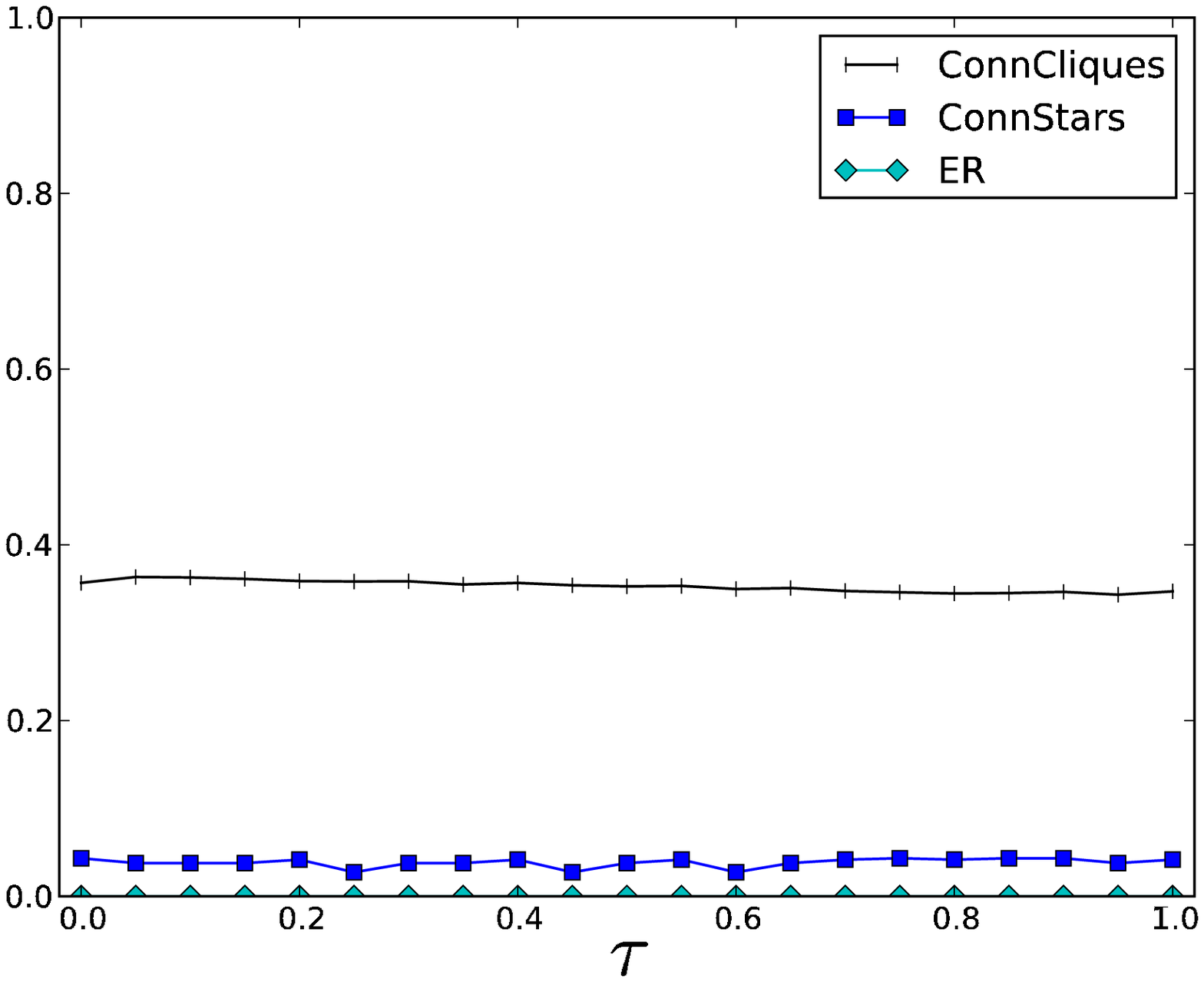}} \subfigure[$r=0.49$]{

\includegraphics[width=0.30\columnwidth]{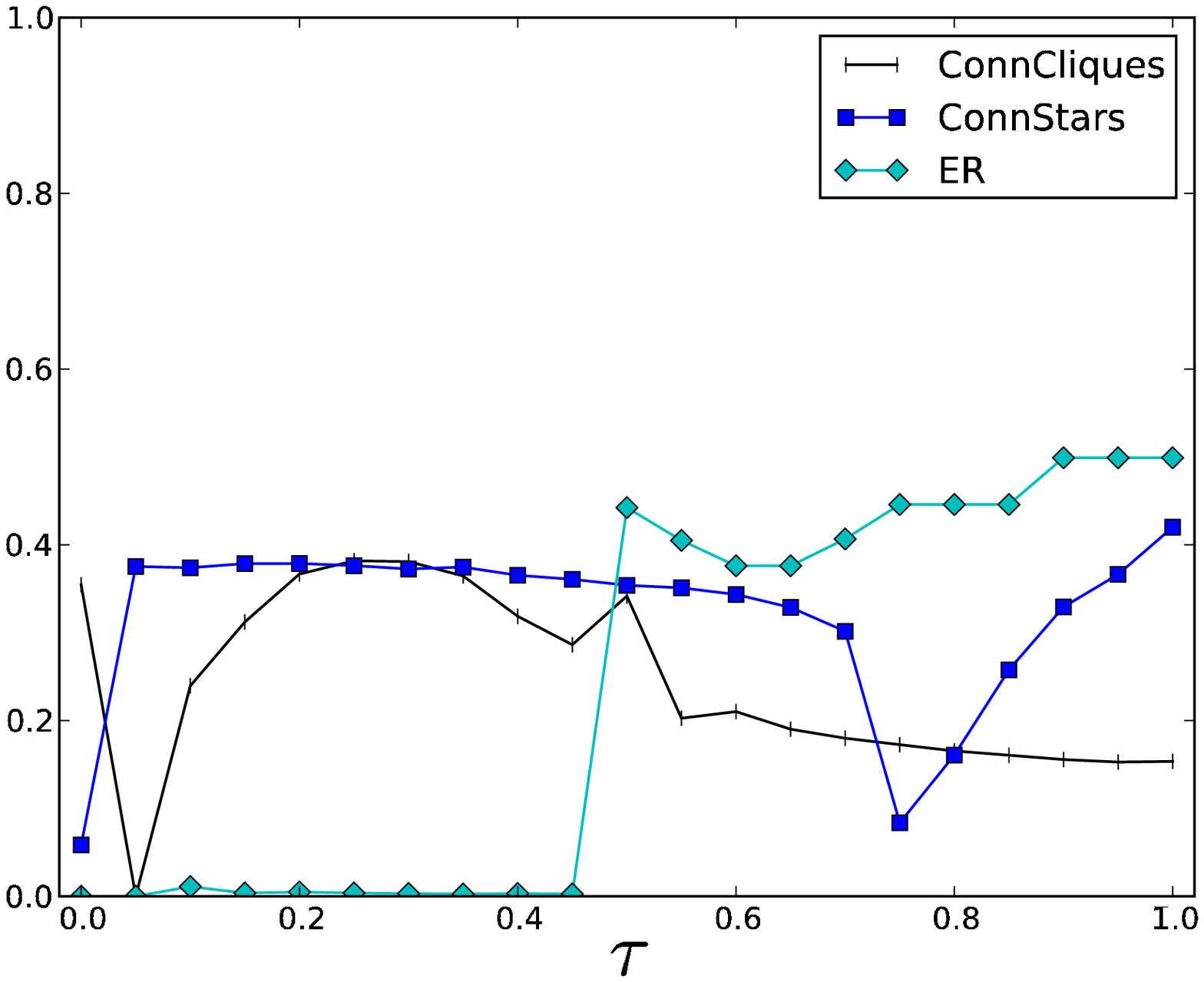}}
\par\end{centering}

\begin{centering}
\subfigure[$r=0.51$]{

\includegraphics[width=0.30\columnwidth]{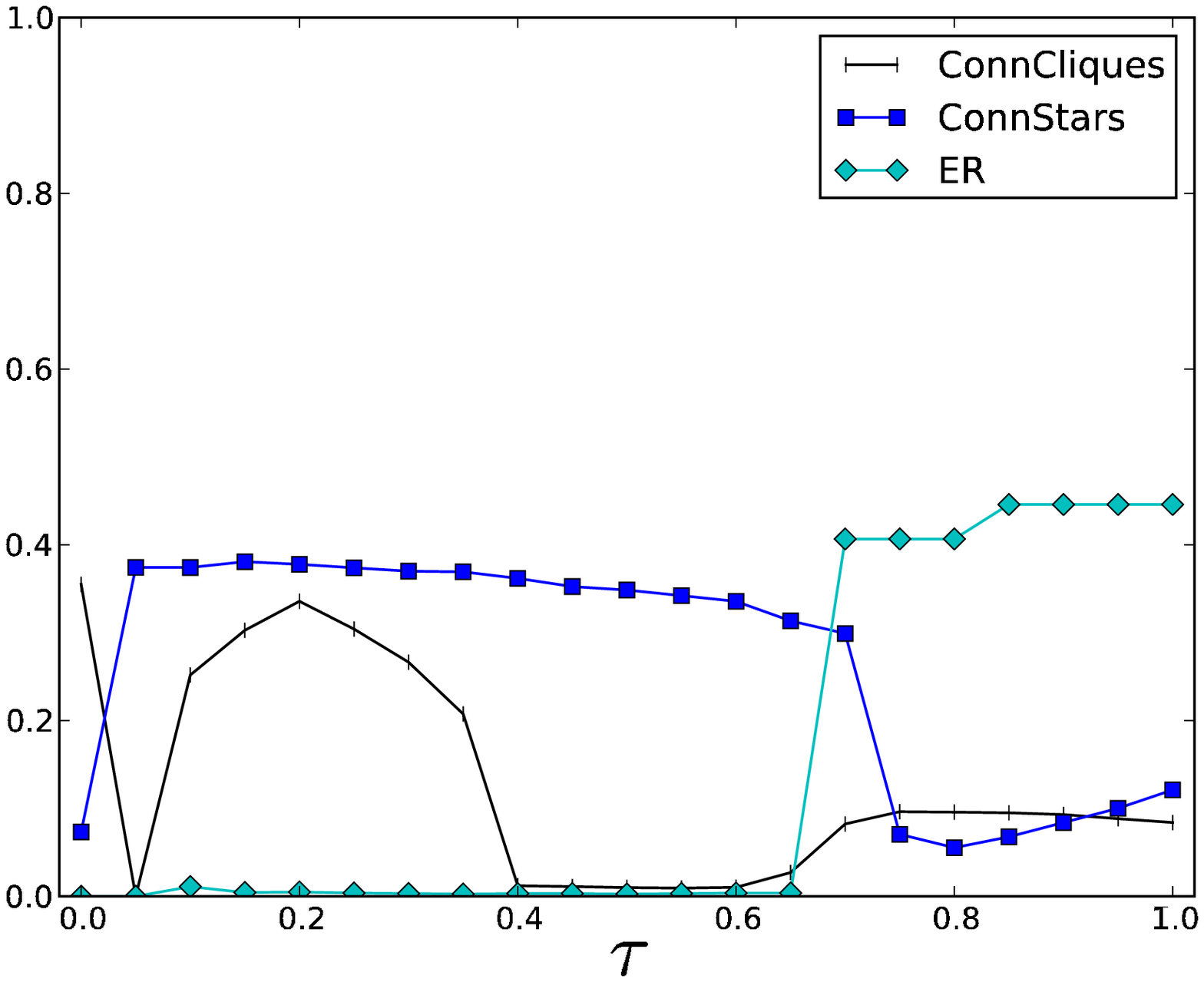}} \subfigure[$r=0.75$]{

\includegraphics[width=0.30\columnwidth]{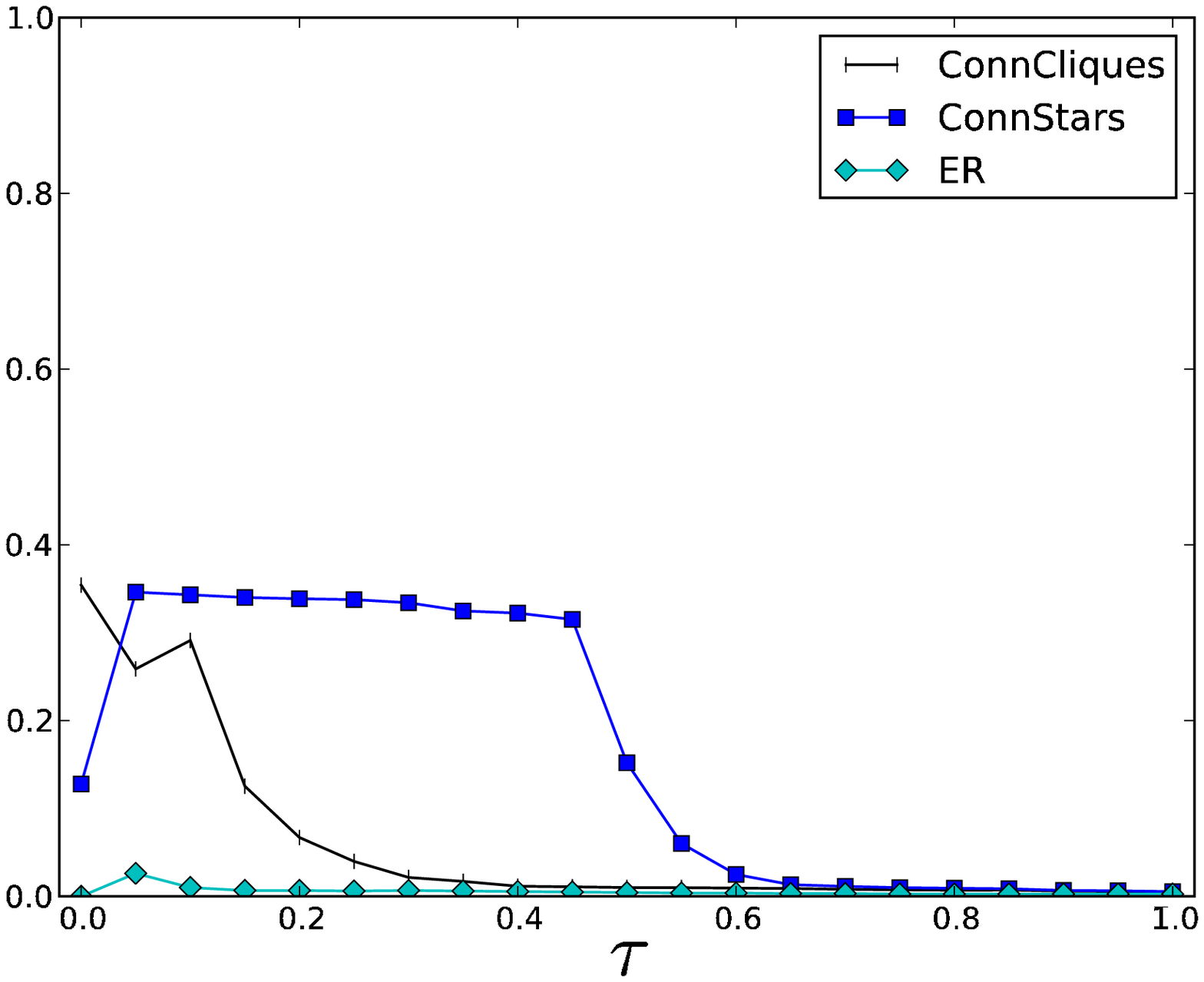}}
\par\end{centering}

\caption{{\bf Standard deviation in connectivity $p$ of the top 5\% of solutions.}}
\label{fig:robustness-opt-p}
\end{figure}

The results are in Figs.\ \ref{fig:robustness-opt-resilience},\ref{fig:robustness-opt-efficiency},\ref{fig:robustness-opt-p},\ref{fig:robustness-opt-k}.
Overall, as one would expect, the properties are more variable near
the transition point $r=0.5$, as compared to $r$ values away from
$r=0.5$. Moreover, variability is high within each design whenever
the design undergoes a phase transition, since multiple different
phases have nearly equal fitness. Another source of variance is found
in designs with two parameters. These are more variable than designs
with a single parameter because the former can reproduce some of the
same graphs with many different parameter settings - the parameters
have non-orthogonal effects. Finally, designs with low fitness for
all configurations (like Cycles) show high variability since all configurations
show uniform low performance.

\section*{Acknowledgments}
This work has benefited from discussions with Aaron Clauset, Michael Genkin, Vadas Gintautas, Shane Henderson, Jason Johnson and Roy Lindelauf, and
two anonymous reviewers.
Part of this work was funded by the Department of Energy at the Los
Alamos National Laboratory (LA-UR 10-01563) under contract DE-AC52-06NA25396 through
the Laboratory Directed Research and Development program, and by the
Defense Threat Reduction Agency.

\bibliography{gutfraind}

\end{document}